# On the Theory of Ends of a Pro-$p$ Group

by Kay Wingberg

**Abstract.** We study the group $H^1(G, \mathbb{F}_p[\![G]\!])$ of ends of a pro-$p$ group $G$ and prove a pro-$p$ analog of Stallings' decomposition theorem.

One of the most important results in the theory of abstract groups is Stallings' decomposition theorem [16]. Let $e(G) = 1 + \dim_{\mathbb{F}_2} H^1(G, \mathbb{F}_2[G])$ be the number of ends of an infinite finitely generated abstract group $G$. Then $e(G) = 1, 2$ or $\infty$, and $e(G) = 2$ if and only if $G$ has a subgroup of finite index which is isomorphic to $\mathbb{Z}$ (Hopf 1943). If $G$ is torsion-free, then the theorem of Stallings (1968) asserts that $e(G) = \infty$ if and only if $G$ is a free product of non-trivial subgroups. Up to now it was impossible to obtain an analog of this theorem for profinite groups. But we will see that for pro-$p$ groups the situation is much better.

If $p$ is a fixed prime number and $G$ a pro-$p$ group, then $\Lambda_G$ denotes the completed group ring $\mathbb{F}_p[\![G]\!]$ of $G$ over $\mathbb{F}_p$ and $I_G$ is its augmentation ideal, i.e. the kernel of the augmentation map $\Lambda_G \twoheadrightarrow \mathbb{F}_p$. Consider the (continuous) cohomology groups $H^i(G, \mathbb{F}_p[\![G]\!])$ and let

$$h^i(G) = \dim_{\mathbb{F}_p} H^i(G, \mathbb{F}_p[\![G]\!]).$$

The number

$$e(G) = 1 - h^0(G) + h^1(G)$$

is called the *number of ends* of $G$. If $G$ is infinite, then $e(G) = 1 + h^1(G)$.

We denote the coinvariants of a left (resp. right) $\Lambda_G$-module $M$ by $M_G = M/I_G M$ (resp. $M_G = M/M I_G$). Considering the right $\Lambda_G$-module $H^1(G, \mathbb{F}_p[\![G]\!])$, let

$$f(G) = \dim_{\mathbb{F}_p} H^1(G, \mathbb{F}_p[\![G]\!])_G.$$

For pro-$p$ groups first results were obtained by Korenev [9]. He proved the following analog of Hopf's theorem:

**Theorem 1:** *Let $G$ be a pro-$p$ group, then*

(i) $e(G) = 0, 1, 2$ *or* $\infty$.
(ii) $e(G) = 0$ *if and only if $G$ is finite.*
(iii) $e(G) = 2$ *if and only if $G$ has an open subgroup isomorphic to $\mathbb{Z}_p$.*

*In particular, if $G$ is infinite, then $h^1(G) = 0, 1$ or $\infty$, and if $G$ is torsion-free, then $h^1(G) = 1$ if and only if $G \cong \mathbb{Z}_p$.*

We call a non-trivial pro-$p$ group $G$ **freely decomposable** if it is the free pro-$p$ product of $s$ non-trivial closed subgroups with $s > 1$; otherwise $G$ is called **freely indecomposable**. By $s(G)$ we denote the number of freely indecomposable factors of $G$ (which is well-defined, see theorem (2.1) below).

We will give an alternative short proof of Korenev's theorem, but our main result will be the following pro-$p$ analog of Stallings' theorem.

**Theorem 2:** *Let $G$ be a finitely generated pro-$p$ group. Then the following two assertions are equivalent:*

(i) *$G$ is freely decomposable.*
(ii) *$I_G$ is decomposable as left $\mathbb{F}_p[\![G]\!]$-module.*

*Furthermore, the following are equivalent:*

(iii) *There exists an open subgroup $H$ of $G$ such that every open subgroup of $H$ is freely decomposable.*
(iv) $\dim_{\mathbb{F}_p} H^1(G, \mathbb{F}_p[\![G]\!]) = \infty$.

*If, in addition, $G$ is torsion-free, then all four statements are equivalent.*

**Corollary 1:** *Let $G$ be a finitely generated torsion-free pro-$p$ group. Then the following assertions hold:*

1. *If $G$ is not free, then $G$ is the free pro-$p$ product of $s(G) = f(G) + 1$ freely indecomposable closed subgroups; furthermore, the $\mathbb{F}_p[\![G]\!]$-module $H^1(G, \mathbb{F}_p[\![G]\!])$ is free of rank $f(G)$. If $G$ is free, then $s(G) = f(G)$ and there is an exact sequence $0 \to \mathbb{F}_p[\![G]\!] \to \mathbb{F}_p[\![G]\!]^{f(G)} \to H^1(G, \mathbb{F}_p[\![G]\!]) \to 0$.*
2. *Let $H$ be an open subgroup of $G$. Then*

$$s(H) = (G : H)(s(G) - 1) + 1.$$

*In particular, $G$ is freely indecomposable if and only if $H$ is, and $I_G$ is indecomposable as left $\mathbb{F}_p[\![G]\!]$-module if and only if $I_H$ is indecomposable.*

We recall the notion of a pro-$p$ duality group [12] (3.4.6): For a discrete $G$-module $A$ and $i \geq 0$, let

$$D_i(G, A) = \varinjlim_{U} H^i(U, A)^{\vee},$$

where $^{\vee}$ denotes the Pontryagin-dual, the direct limit is taken over all open subgroups $U$ of $G$ and the transition maps are the duals of the corestriction maps. $D_i(G, A)$ is a discrete $G$-module in a natural way. Assume that the cohomological dimension $\operatorname{cd}_p G = n$ is finite. Then the $G$-module

$$D(G) = \varinjlim_{\nu \in \mathbb{N}} D_n(G, \mathbb{Z}/p^{\nu}\mathbb{Z})$$

is called the dualizing module of $G$ and we have in a natural way the trace map

$$tr : H^n(G, D(G)) \longrightarrow \mathbb{Q}_p/\mathbb{Z}_p.$$

The pro-$p$ group $G$ is called a **duality group of dimension $n$** if for all $i \in \mathbb{Z}$ and all finite $p$-primary $G$-modules $A$, the cup-product and the trace map

$$H^i(G, \operatorname{Hom}(A, D(G))) \times H^{n-i}(G, A) \xrightarrow{\cup} H^n(G, D(G)) \xrightarrow{tr} \mathbb{Q}_p/\mathbb{Z}_p$$

yield an isomorphism $H^i(G, \operatorname{Hom}(A, D(G))) \cong H^{n-i}(G, A)^{\vee}$.

**Corollary 2:** *Let $G$ be a finitely generated pro-$p$ group of cohomological dimension $\operatorname{cd}_p G \leq 2$. Then $G$ is the free pro-$p$ product of finitely many duality groups and the following assertions are equivalent:*

(a) *$G$ is a duality group of dimension 2.*

(b) *$G$ is freely indecomposable.*

(c) *$I_G$ is indecomposable as left $\mathbb{F}_p[\![G]\!]$-module.*

*In particular, if $G$ is a 2-generator group with $\operatorname{cd}_p G = 2$, then $G$ is a duality group.*

**Remarks:**
1. If $G$ is a finitely generated FAB pro-$p$ group, i.e. the abelianization $U^{ab}$ of every open subgroup $U$ of $G$ is finite, then, by the principal ideal theorem, $H^1(G, \mathbb{F}_p[\![G]\!]) = 0$, see (1.3)(iii) below. One easily sees directly that $G$ is freely indecomposable: suppose that $G = H_1 * H_2$ is a non-trivial decomposition and let $U$ be the kernel of the map $G \twoheadrightarrow \mathbb{Z}/p\mathbb{Z}$, where each $H_i$ surjects onto $\mathbb{Z}/p\mathbb{Z}$. By the pro-$p$ analog of Kurosh' subgroup theorem $U$ has a free factor of rank $r = p - 1$, hence $U^{ab}$ surjects onto $(\mathbb{Z}_p)^r$, a contradiction.

2. In number theory an example of a class of pro-$p$ groups satisfying the first assertion of corollary 2 was known before: let $G_S = Gal(k_S(p)|k)$ be the Galois group of the maximal $p$-extension of a number field $k$ unramified outside the set $S$ of primes of $k$. If $p \neq 2$, $k$ contains the group of $p$-th roots of unity and $S$ is finite and contains all primes above $p$ and all archimedean primes, then $\mathrm{cd}_p G \leq 2$ and $G_S$ is a free pro-$p$ product of duality groups, i.e. it is itself a duality group or a free pro-$p$ product of decomposition groups and a free pro-$p$ group, see [12] (10.9.8).

3. In an appendix we collect corresponding results in the case of abstract groups.

4. After finishing this paper we learned that Th. Weigel and P. A. Zalesskiĭ also proved an analogue of Stallings' decomposition theorem (for arbitrary finitely generated pro-$p$ groups) using a slightly different defininition for the group of ends, see [19]. Finally I would like to thank Thomas Weigel and John MacQuarrie for pointing out some errors in an earlier version of this paper.

# 1 Decomposition of $I_G$

Let $p$ be a fixed prime number, $G$ a pro-$p$ group and $k$ a finite field of characteristic $p$. Then $\Lambda_G(k) = k[\![G]\!]$ denotes the completed group ring of $G$ over $k$ with augmentation ideal $I_G(k)$; $\Lambda_G(k)$ is a local ring and indecomposable as $\Lambda_G(k)$-module. Observe that a projective $\Lambda_G(k)$-module is free. If $k = \mathbb{F}_p$, we write for short $\Lambda_G = \Lambda_G(\mathbb{F}_p)$ and $I_G = I_G(\mathbb{F}_p)$.

For the proof of the following result we refer to [18] prop.(2.1).

**Proposition 1.1** *Let $G$ be a pro-$p$ group.*

(i) *Let $M$ be a non-zero indecomposable finitely generated left $\Lambda_G(k)$-module. Then $\mathrm{End}_{\Lambda_G(k)}(M)$ is a local ring.*

(ii) *The Krull-Schmidt-Azumaya theorem holds, i.e. if $M$ is a finitely generated left $\Lambda_G(k)$-module, then $M$ is expressible as a finite direct sum of indecomposable left submodules. Further, if*

$$M = \bigoplus_{i=1}^{r} M_i = \bigoplus_{j=1}^{s} N_j,$$

*are two such sums, then $r = s$ and, by re-ordering, we have $M_i \cong N_i$, $i = 1, \ldots, r$.*

We introduce a notation which is used frequently in representation theory: write
$$N | M$$

to mean that the $\Lambda_G(k)$-module $N$ is isomorphic to a direct summand of the $\Lambda_G(k)$-module $M$.

**Corollary 1.2** *Let $G$ be a pro-$p$ group, $N_1$, $N_2$ and $M$ left $\Lambda_G(k)$-modules where $M$ is non-zero, finitely generated and indecomposable.*
*Then $M|(N_1 \oplus N_2)$ implies $M|N_1$ or $M|N_2$.*

**Proof:** This is a formal consequence of the fact that $\mathrm{End}_{\Lambda_G(k)}(M)$ is a local ring, see for example [10] (4.5). □

Let $G^* = G^p[G, G]$ be the Frattini-subgroup of $G$, then

$$d(G) = \dim_{\mathbb{F}_p} G/G^* = \dim_{\mathbb{F}_p} H^1(G, \mathbb{Z}/p\mathbb{Z})$$

is the minimal number of generators of $G$. If $V \subseteq U$ are open normal subgroups of $G$, then

$$\mathrm{Ver}_U^V \colon U/U^* \longrightarrow (V/V^*)^U$$

is the transfer map.

Let $k$ be finite field and let $M$ be a compact left $\Lambda_G(k)$-module and $N$ a compact $\Lambda_G(k)$-bimodule. Then the group $\mathrm{Hom}_{\Lambda_G(k)}(M, N)$ of continuous $G$-homomorphisms with the compact-open topology becomes a compact right $\Lambda_G(k)$-module by $(\varphi g)(m) = \varphi(m)g$, $\varphi \in \mathrm{Hom}_{\Lambda_G(k)}(M, N), g \in G, m \in M$. The continuous cohomology group $H^1(G, M)$ has a natural structure as a right $\Lambda_G(k)$-module which is induced by the action on the group of 1-cocycles: Let $a\colon G \to M$ be a 1-cocycle and $g \in G$, then $(a \cdot g)(x) = a(x)g$ for $x \in G$.

**Lemma 1.3** *Let $G$ be a pro-$p$ group.*

(i) *There is a natural $\Lambda_G$-isomorphism*

$$H^1(G, \mathbb{F}_p[\![G]\!]) \cong \varprojlim_U H^1(G, \mathbb{F}_p[G/U]),$$

*where $U$ runs through the normal open subgroups of $G$.*

(ii) *Let $H$ be an open subgroup of $G$. Then there is a natural $\Lambda_H$-isomorphism*

$$H^1(G, \mathbb{F}_p[\![G]\!]) \xrightarrow[\sim]{sh} H^1(H, \mathbb{F}_p[\![H]\!]).$$

(iii) *If $G$ is finitely generated, then there are natural $\Lambda_G$-isomorphisms*

$$H^1(G, \mathbb{F}_p[\![G]\!]) \cong (\varinjlim_{U, \mathrm{Ver}_G^U} U/U^*)^\vee \cong D_1(G, \mathbb{F}_p)^\vee,$$

*where the limit is taken over the transfer maps.*

**Proof:** Using [12] (2.7.6) and Shapiro's lemma [12] (1.6.4), we have

$$H^1(G, \mathbb{F}_p[\![G]\!]) \cong \varprojlim_U H^1(G, \mathbb{F}_p[G/U]) \cong \varprojlim_U H^1(U, \mathbb{F}_p),$$

$$H^1(G, \mathbb{F}_p[\![G]\!]) = H^1(G, Ind_H^G \mathbb{F}_p[\![H]\!]) \stackrel{sh}{\Rightarrow} H^1(H, \mathbb{F}_p[\![H]\!]).$$

If $G$ is finitely generated, then the groups $H^1(G, \mathbb{F}_p[G/U])$ are finite and the $\Lambda_G$-module $H^1(G, \mathbb{F}_p[\![G]\!])$ is compact. Pontryagin-duality gives

$$H^1(G, \mathbb{F}_p[\![G]\!]) \cong (\varinjlim_U H^1(U, \mathbb{F}_p)^\vee)^\vee \cong (\varinjlim_U U/U^*)^\vee \cong D_1(G, \mathbb{F}_p)^\vee.$$

$\square$

Let $G$ be an infinite pro-$p$ group. Then $\operatorname{Hom}_{\Lambda_G}(\mathbb{F}_p, \Lambda_G) = 0$. Since

$$\operatorname{Ext}^1_{\Lambda_G}(\mathbb{F}_p, \Lambda_G) = H^1(G, \Lambda_G) \quad \text{and} \quad \operatorname{Hom}_{\Lambda_G}(\Lambda_G, \Lambda_G) = \Lambda_G,$$

see [12] (5.2.14), the exact sequence

$$0 \longrightarrow I_G \longrightarrow \Lambda_G \longrightarrow \mathbb{F}_p \longrightarrow 0$$

yields the exact sequence

(∗) $\qquad 0 \longrightarrow \Lambda_G \xrightarrow{\phi_G} \operatorname{Hom}_{\Lambda_G}(I_G, \Lambda_G) \longrightarrow H^1(G, \Lambda_G) \longrightarrow 0,$

where $\phi_G(\lambda) : x \mapsto x\lambda$.

**Proposition 1.4** *Let $G$ be a pro-$p$ group.*

(i) *If $I_G$ is decomposable, then $G$ is infinite and $h^1(G) \neq 0$.*

(ii) *If $h^1(G) \neq 0$, then there exists an open subgroup $H$ of $G$ such that $I_H$ has a direct summand isomorphic to $\Lambda_H$.*

(iii) *If $I_G \cong \Lambda_G$, then $G \cong \mathbb{Z}_p$ and $h^1(G) = 1$.*

(iv) *If $0 < h^1(G) < \infty$, then $h^1(G) = 1$ and there exists an open subgroup $H$ of $G$ isomorphic to $\mathbb{Z}_p$.*

**Proof:** (i) Let $I_G = M_1 \oplus M_2$ be a non-trivial decomposition of $I_G$. If $G$ is finite, then $(M_1)^G \oplus (M_2)^G = (I_G)^G$ and $\dim_{\mathbb{F}_p}(I_G)^G = 1$ implies that $M_1 = 0$ or $M_2 = 0$. Thus $G$ has to be infinite. We consider the exact sequence

$$0 \longrightarrow \Lambda_G \longrightarrow \operatorname{Hom}_{\Lambda_G}(M_1, \Lambda_G) \oplus \operatorname{Hom}_{\Lambda_G}(M_2, \Lambda_G) \longrightarrow H^1(G, \Lambda_G) \longrightarrow 0.$$

Since the groups $\operatorname{Hom}_{\Lambda_G}(M_j, \Lambda_G)$, $j = 1, 2$, are non-zero and $\Lambda_G$ is indecomposable, assertion (i) follows.

(ii) Since $H^1(G, \Lambda_G) = \varprojlim_U H^1(U, \mathbb{F}_p) \neq 0$, see (1.3)(i), and since the diagram

$$\begin{array}{ccc} H^1(G, \Lambda_G) & \longrightarrow & H^1(G, \mathbb{F}_p) \\ {\scriptstyle sh} \downarrow \sim & & \uparrow {\scriptstyle cor} \\ H^1(U, \Lambda_U) & \longrightarrow & H^1(U, \mathbb{F}_p) \end{array}$$

commutes, see [12] (1.6.5), there exists an open subgroup $H$ of $G$ such that the canonical map $H^1(H, \mathbb{F}_p[\![H]\!]) \to H^1(H, \mathbb{F}_p)$ induce by the augmentation map $\mathbb{F}_p[\![H]\!] \to \mathbb{F}_p$ is not zero. Now the commutative and exact diagram

$$\begin{array}{ccc} \operatorname{Hom}_{\Lambda_H}(I_H, \mathbb{F}_p) & \xrightarrow{\sim} & H^1(H, \mathbb{F}_p) \\ \uparrow & & \uparrow \\ 0 \longrightarrow \Lambda_H \longrightarrow \operatorname{Hom}_{\Lambda_H}(I_H, \Lambda_H) & \longrightarrow & H^1(H, \Lambda_H) \longrightarrow 0 \end{array}$$

shows that there exists a $\Lambda_H$-homomorphism $I_H \to \Lambda_H$ which is surjective as its image is not contained in $I_H$. It follows that $I_H \cong \Lambda_H \oplus M$ with some $\Lambda_H$-module $M$ (possibly trivial).

(iii) The exact sequence $0 \to \Lambda_G \to \Lambda_G \to \mathbb{F}_p \to 0$ shows that the projective dimension of $\mathbb{F}_p$ as $\Lambda_G$-module is equal to 1, hence $\operatorname{cd}_p G = 1$, see [12] (5.2.13). Since $I_G \cong (\Lambda_G)^r$ for a free pro-$p$ group $G$ of rank $r$, [12] (5.6.3), (5.6.4), we obtain $G \cong \mathbb{Z}_p$, and so $h^1(G) = 1$.

(iv) Since $0 < h^1(G) < \infty$, $G$ has to be infinite by Shapiro's lemma. Using (ii), we see that there exists an open subgroup $H$ of $G$ such that $I_H$ has a direct summand isomorphic to $\Lambda_H$, i.e. $I_H \cong \Lambda_H \oplus M$ for some $\Lambda_H$-module $M$. We get the commutative and exact diagram

$$\begin{array}{c} \operatorname{Hom}_{\Lambda_H}(M, \Lambda_H) \\ \cap \\ \downarrow \\ 0 \longrightarrow \Lambda_H \xrightarrow{\phi_H} \operatorname{Hom}_{\Lambda_H}(I_H, \Lambda_H) \longrightarrow H^1(H, \Lambda_H) \longrightarrow 0 \\ \| \qquad\qquad \downarrow \\ \Lambda_H \xrightarrow{\varepsilon} \operatorname{Hom}_{\Lambda_H}(\Lambda_H, \Lambda_H). \end{array}$$

The image of $\phi_H$ consists of homomorphisms $I_H \to \Lambda_H$ given by right-multiplication with an element $\lambda \in \Lambda_H$. Since the $\Lambda_H$-annulator of $\Lambda_H$ is zero, the map $\varepsilon$ is injective. It follows that $\operatorname{Hom}_{\Lambda_H}(M, \Lambda_H)$ injects into $H^1(H, \Lambda_H) \cong H^1(G, \Lambda_G)$. Hence this module is finite and so fixed by an open subgroup $H_0$ of $H$. Therefore the map $\iota \in \operatorname{Hom}_{\Lambda_H}(M, \Lambda_H)$, $\iota: M \hookrightarrow I_H \hookrightarrow \Lambda_H$, has an image in $(\Lambda_H)^{H_0} = 0$, hence $M = 0$. Thus $I_H \cong \Lambda_H$, and so $H \cong \mathbb{Z}_p$ and $1 = h^1(H) = h^1(G)$ by (iii). This proves (iv). $\square$

**Proof of Theorem 1:** Obviously $h^0(G) \leq 1$ and $h^0(G) = 1$ if and only if $G$ is finite. Using Shapiro's lemma, we obtain assertion (ii).

Assume that $0 < h^1(G) < \infty$. Then $G$ has to be infinite by Shapiro's lemma and $e(G) = h^1(G) + 1$. From (1.4)(iv) it follows that $h_1(G) = 1$ and the exists an open subgroup $H$ of $G$ isomorphic to $\mathbb{Z}_p$, i.e. the assertions (i) and (iii) hold.

If $G$ is torsion-free, then a theorem of Serre, see [15], implies that $\mathrm{cd}_p G = 1$, and so $G \cong \mathbb{Z}_p$. □

**Lemma 1.5** *Let $G$ and $H$ be pro-$p$ groups and $k$ a finite field of characteristic $p$. Let $M$ be a left compact $\Lambda_H(k)$-module, $L$ a compact $(\Lambda_G(k), \Lambda_H(k))$-bimodule and $N$ a left $\Lambda_G(k)$-module. Assume either that $N$ is compact and $M$ and $L$ are finitely generated or that $N$ is discrete.*

(i) *There is a canonical isomorphism*

$$S : \mathrm{Hom}_{\Lambda_H(k)}(M, \mathrm{Hom}_{\Lambda_G(k)}(L, N)) \xrightarrow{\sim} \mathrm{Hom}_{\Lambda_G(k)}(L \hat{\otimes}_{\Lambda_H(k)} M, N),$$

*such that $(S\varphi)(l \hat{\otimes} m) = \varphi(m)(l)$. This morphism induces a natural equivalence of functors.*

(ii) *Let $H$ be a closed subgroup of $G$. Considering $\Lambda_G(k)$ as $(\Lambda_G(k), \Lambda_H(k))$-bimodule, there are canonical isomorphisms*

$$\mathrm{Ext}^i_{\Lambda_H(k)}(M, Res_H N) \xrightarrow{\sim} \mathrm{Ext}^i_{\Lambda_G(k)}(\Lambda_G(k) \hat{\otimes}_{\Lambda_H(k)} M, N), \ i \geq 0,$$

*where $Res_H N$ denotes the $\Lambda_H(k)$-module $N$ obtained by restriction of scalars.*

**Proof:** If $A$ and $B = \varprojlim_i B_i$ are compact $\Lambda_G(k)$-modules ($B_i$ finite $\Lambda_G(k)$-modules), where $A$ is finitely generated, then

$$\mathrm{Hom}_{\Lambda_G(k)}(A, B) = \varprojlim_i \mathrm{Hom}_{\Lambda_G(k)}(A, B_i)$$

is a compact $\Lambda_G(k)$-module.

Now assertion (i) is the topological analog of [4] chap.II (5.2) or [6] (2.19): see [3](2.4), where N is discrete, and take the projective limit if $N$ is compact.

In order to prove (ii), take $L = \Lambda_G(k)$ and observe that $\mathrm{Hom}_{\Lambda_G(k)}(\Lambda_G(k), N)$ and $Res_H N$ are isomorphic as left $\Lambda_H(k)$-modules. Then we obtain from (i) the assertion for $i = 0$. Since $\Lambda_G(k)$ is $\Lambda_H(k)$-projective, the functor $\Lambda_G(k) \hat{\otimes}_{\Lambda_H(k)} -$ is exact, and so, taking a $\Lambda_H(k)$-projective resolution of $M$, we obtain the desired isomorphisms for all $i \geq 0$, see also [4] VI. Proposition (4.1.3). □

Let $G$ be a pro-$p$ group and let $H$ be a closed subgroup. If $I$ is a left ideal of $\Lambda_H$, then $\Lambda_G \hat{\otimes}_{\Lambda_H} I \cong \Lambda_G I$, see [17] (4.3): Applying the exact functor $\Lambda_G \hat{\otimes}_{\Lambda_H} -$ to the exact sequence $0 \to I \to \Lambda_H$ gives $0 \to \Lambda_G \hat{\otimes}_{\Lambda_H} I \to \Lambda_G$ with image $\Lambda_G I$. In particular,
$$J_H := \Lambda_G \hat{\otimes}_{\Lambda_H} I_H = \mathrm{Ind}_H^G I_H$$
is the left ideal of $\Lambda_G$ generated by $I_H$. The Frobenius reciprocity (1.5)(ii) gives the natural isomorphisms
$$\mathrm{Ext}^i_{\Lambda_H}(I_H, \mathrm{Res}_H N) \cong \mathrm{Ext}^i_{\Lambda_G}(J_H, N)$$
for $i \geq 0$, where $N$ is a discrete left $\Lambda_G$-module or $N$ is compact and $H$ is finitely generated.

The following result is a pro-$p$ analog of a theorem for abstract groups, see [5] theorem (4.7).

**Proposition 1.6** *Let $G$ be a pro-$p$ group, and let $H_j$, $1 \leq j \leq n$, be closed subgroups. Then the following assertions are equivalent*

(i) $\qquad\qquad\qquad\qquad G = H_1 * \cdots * H_n.$

(ii) $\qquad\qquad\qquad\qquad I_G = J_{H_1} \oplus \cdots \oplus J_{H_n}.$

**Proof:** By [12] (4.1.5) we know that (i) is equivalent to

(iii) $\qquad\qquad H^i(G, \mathbb{F}_p) \xrightarrow[\sim]{res} \bigoplus_{j=1}^n H^i(H_j, \mathbb{F}_p), \ i = 1, 2.$

From the exact sequence $0 \to I_G \to \Lambda_G \to \mathbb{F}_p \to 0$ we obtain the isomorphisms
$$\mathrm{Ext}^i_{\Lambda_G}(I_G, \mathbb{F}_p) \xrightarrow{\sim} \mathrm{Ext}^{i+1}_{\Lambda_G}(\mathbb{F}_p, \mathbb{F}_p) = H^{i+1}(G, \mathbb{F}_p), \ i = 0, 1.$$
Since $\mathrm{Ext}^i_{\Lambda_G}(J_{H_j}, \mathbb{F}_p) \cong \mathrm{Ext}^i_{\Lambda_{H_j}}(I_{H_j}, \mathbb{F}_p)$, assertion (iii) is equivalent to

(iv) $\qquad\qquad \mathrm{Ext}^i_{\Lambda_G}(I_G, \mathbb{F}_p) \xrightarrow[\sim]{can} \bigoplus_{j=1}^n \mathrm{Ext}^i_{\Lambda_G}(J_{H_j}, \mathbb{F}_p), \ i = 0, 1.$

Therefore (ii) implies (i). Conversely, from the exact sequence
$$0 \longrightarrow M \longrightarrow \bigoplus_{j=1}^n J_{H_j} \longrightarrow I_G \longrightarrow 0,$$
where $M$ is defined as kernel of the natural surjection on the right, we obtain the exact sequence
$$0 \to \mathrm{Hom}_{\Lambda_G}(I_G, \mathbb{F}_p) \to \bigoplus_{j=1}^n \mathrm{Hom}_{\Lambda_G}(J_{H_j}, \mathbb{F}_p) \to \mathrm{Hom}_{\Lambda_G}(M, \mathbb{F}_p)$$
$$\to \mathrm{Ext}^1_{\Lambda_G}(I_G, \mathbb{F}_p) \to \bigoplus_{j=1}^n \mathrm{Ext}^1_{\Lambda_G}(J_{H_j}, \mathbb{F}_p).$$

Using (iv), we obtain $\text{Hom}_{\Lambda_G}(M, \mathbb{F}_p) = 0$, hence $M = 0$. This completes the proof of the proposition. □

The following corollary can also be obtained from theorem (3.4) in [13].

**Corollary 1.7** *Let $G$ be a finitely generated pro-p duality group of dimension $n \geq 2$. Then $G$ is freely indecomposable.*

**Proof:** By [12] (3.4.6) we have $D_1(G, \mathbb{F}_p) = 0$, if $G$ is a duality group of dimension $n \geq 2$. Using lemma (1.3)(iii), it follows that $h^1(G) = 0$. Hence proposition (1.6) and proposition (1.4)(i) gives the result. □

In order to establish a decomposition of the form $I_G = J_{H_1} \oplus \cdots \oplus J_{H_n}$ we will need the following lemmata.

**Lemma 1.8** *Let $U$ be a closed subgroup of the pro-p group $G$ and let $k$ be a finite field of characteristic $p$. Then*

$$Res_U I_G(k) = I_U(k) \oplus P,$$

*where $P$ is a free $\Lambda_U(k)$-module. If $U$ is open and normal in $G$ of index $d$ and if $\{\eta_1, \ldots, \eta_{d-1}\} \subseteq I_G(k)$ is an arbitrary pre-image of a $k$-basis of $I_{G/U}(k)$ under the canonical surjection $I_G(k) \twoheadrightarrow I_{G/U}(k)$, then*

$$Res_U I_G(k) = I_U(k) \oplus \bigoplus_{i=1}^{d-1} \Lambda_U(k)\, \eta_i.$$

**Proof:** Consider the commutative and exact diagram

$$\begin{array}{ccccccccc}
0 & \longrightarrow & Res_U I_G(k) & \longrightarrow & Res_U \Lambda_G(k) & \longrightarrow & k & \longrightarrow & 0 \\
& & \uparrow & & \uparrow \varphi & & \| & & \\
0 & \longrightarrow & I_U(k) & \longrightarrow & \Lambda_U(k) & \longrightarrow & k & \longrightarrow & 0,
\end{array}$$

where $\varphi(\Lambda_U(k))$ is some free summand (of rank 1) of the free $\Lambda_U(k)$-module $Res_U \Lambda_G(k)$. It follows that $Res_U I_G(k)$ surjects onto the free $\Lambda_U(k)$-module $P = Res_U \Lambda_G(k)/\varphi(\Lambda_U(k))$, showing the first statement. Since $Res_U \Lambda_G(k) = \bigoplus_{i=0}^{d-1} \Lambda_U(k)\, \eta_i$, $\eta_0 = 1$, the second is also obvious. □

Recall that a $\Lambda_G(k)$-module $M$ is called $G$-$H$-projective, where $H$ is a closed subgroup of the profinite group $G$, iff ever a $\Lambda_G(k)$-epimorphism $B \twoheadrightarrow M$ splits as a $\Lambda_H(k)$-module homomorphism, then it splits as a $\Lambda_G(k)$-module homomorphism, see [10] (3.1),(3.2).

**Lemma 1.9** *Let $G$ be a finitely generated pro-$p$ group and $H$ a closed subgroup of $G$. Let $M$ be a finitely generated $\Lambda_G$-module and $k|\mathbb{F}_p$ a finite extension of fields. Then*

(i) $\qquad\qquad M$ *is $G$-$H$-projective* $\Leftrightarrow k\hat\otimes_{\mathbb{F}_p} M$ *is $G$-$H$-projective.*

(ii) *Assume in addition that $M|I_G$ and $M$ is indecomposable and not isomorphic to $\Lambda_G$. Then*
$$M|Ind_H^G I_H \Leftrightarrow k\hat\otimes_{\mathbb{F}_p} M | Ind_H^G I_H(k).$$

**Proof:** (i) If $A$ is a $\Lambda_G$-module, then we put $\bar A = k\hat\otimes_{\mathbb{F}_p} A$. Let $\bar M$ be $G$-$H$-projective and let

(∗) $\qquad\qquad 0 \longrightarrow A \longrightarrow B \longrightarrow M \longrightarrow 0$

be an exact sequence of $\Lambda_G$-modules which is $H$-split. Applying the exact functor $k\hat\otimes_{\mathbb{F}_p} -$, we obtain the exact sequence $0 \to \bar A \to \bar B \to \bar M \to 0$ of $\Lambda_G(k)$-modules which is $H$-split, hence $G$-split. Therefore we obtain for every open normal subgroup $N$ of $G$ a split exact sequence
$$0 \longrightarrow \bar A_N \longrightarrow \bar B_N \longrightarrow \bar M_N \longrightarrow 0$$
of $k[G/N]$-modules. If follows that the sequences $0 \to A_N \to B_N \to M_N \to 0$ are exact and these sequences split since
$$k \otimes_{\mathbb{F}_p} \operatorname{Ext}_{\mathbb{F}_p[G/H]}(M_N, A_N) \xrightarrow{\sim} \operatorname{Ext}_{k[G/H]}(\bar M_N, \bar A_N),$$
see [6] (8.16). Thus the sequence (∗) is $G$-split.

Conversely, if $M$ is $G$-$H$-projective, then $M|Ind_H^G Res_H M$, [10] (3.2), and so $k\hat\otimes_{\mathbb{F}_p} M | Ind_H^G Res_H(k\hat\otimes_{\mathbb{F}_p} M)$. Again by [10] (3.2), it follows that $k\hat\otimes_{\mathbb{F}_p} M$ is $G$-$H$-projective.

(ii) In order to prove the non-trivial implication, let $k\hat\otimes_{\mathbb{F}_p} M | Ind_H^G I_H(k)$. Then $k\hat\otimes_{\mathbb{F}_p} M$ is $G$-$H$-projective, and so $M$ is $G$-$H$-projective. Using (1.8), we have $Res_H M | I_H \oplus P_0$, where $P_0$ is a free $\Lambda_H$-module, and, using [10] (3.2), we obtain
$$M | Ind_H^G Res_H M | Ind_H^G I_H \oplus P,$$
where $P$ is a free $\Lambda_G$-module. Since $M$ is indecomposable and not free, it follows that $M | Ind_H^G I_H$, see (1.2). $\square$

**Lemma 1.10** *Let $G$ be a finitely generated pro-$p$ group. Assume that there is a decomposition*
$$I_G = P \oplus M$$
*of left $\Lambda_G$-modules where $P$ is free of rank $r$. Then there exists a free pro-$p$ subgroup $F$ of $G$ of rank $r$ such that*
$$I_G = Ind_F^G I_F \oplus M.$$

**Proof:** Let $n = \dim_{\mathbb{F}_p} G/G^*$ and

$$G/G^* = <\overline{x}_1, \ldots, \overline{x}_r> \oplus <\overline{x}_{r+1}, \ldots, \overline{x}_n>$$

be the decomposition of $G/G^*$ corresponding to

$$G/G^* \xrightarrow{can} I_G/I_G^2 = P_G \oplus M_G,$$

where $can$ is induced by the map $\sigma \mapsto \sigma - 1$, and so $P_G = <\overline{x}_i - 1, i = 1, \ldots, r>$. Let $x_1 \ldots, x_n \in G$ be arbitrary pre-images of the $\overline{x}_i$'s, and let $F = <x_1, \ldots, x_r>$. Let $F_n$ be a free pro-$p$ group with basis $\{y_1, \ldots, y_n\}$ and let

$$1 \longrightarrow R \longrightarrow F_n \xrightarrow{\pi} G \longrightarrow 1, \ \pi(y_i) = x_i, \ i = 1, \ldots, n,$$

a be minimal representation of $G$. If $F_r$ denotes the free subgroup of $F_n$ generated by $\{y_1, \ldots, y_r\}$, then $F_r$ is mapped onto $F$. Since

$$I_{F_n} = \bigoplus_{i=1}^n \Lambda_{F_n}(y_i - 1) \ \text{ and } \ I_{F_n}/I_R I_{F_n} = \bigoplus_{i=1}^n \Lambda_G(y_i - 1),$$

see [12] (5.6.4), (5.6.6), we obtain a commutative diagram

$$\begin{array}{c}
\bigoplus_{i=1}^n \Lambda_G(y_i - 1) \xrightarrow{\tilde{\pi}} I_G \xrightarrow{\xi} P \\
\iota \uparrow \qquad \qquad \uparrow \quad \nearrow \psi \\
\bigoplus_{i=1}^r \Lambda_G(y_i - 1) \xrightarrow{\tilde{\pi}} X,
\end{array}$$

with $\varphi$ across the top.

where $\tilde{\pi}$ is induced by $\pi$, $X$ is the image of $(Ind_{F_r}^{F_n} I_{F_r})_R = \bigoplus_{i=1}^r \Lambda_G(y_i - 1)$ under $\tilde{\pi}$, $\xi$ is the surjection of $I_G$ onto the free factor given by assumption and $\varphi$ resp. $\psi$ its composition with $\tilde{\pi}$ resp. the inclusion. The map $\iota\varphi = \psi\tilde{\pi}$ is surjective, and so $\psi$ and $\tilde{\pi}$ restricted to $\bigoplus_{i=1}^r \Lambda_G(y_i - 1)$ are bijective. It follows that $X = Ind_F^G I_F$ is a free $\Lambda_G$-module of rank $r$, and so $F$ is free of rank $r$, and $I_G = X \oplus M$. □

**Lemma 1.11** *Let $G$ be a finitely generated pro-p group and $k$ a finite field of characteristic $p$. Let*

$$I_G(k) = M \oplus R$$

*be a non-trivial decomposition of $I_G(k)$ into left $\Lambda_G(k)$-modules, where $M$ is indecomposable and not free. Then there exists a proper open normal subgroup $E$ of $G$ such that*

$$M | Ind_E^G I_E(k) \quad \text{and} \quad M + I_G(k)^2 = I_E(k) + I_G(k)^2.$$

**Proof:** Let $E$ be an open normal subgroup of $G$ containing $G^*$ such that $(E/G^*)\hat\otimes_{\mathbb{F}_p} k$ corresponds to $M_G$ in the decomposition

$$G/G^* \hat\otimes_{\mathbb{F}_p} k \xrightarrow{can} I_G(k)/I_G(k)^2 = M_G \oplus R_G.$$

Then $d = (G : E) > 1$ since the decomposition of $I_G(k)$ is non-trivial and

$$M_G \cong M + I_G(k)^2/I_G(k)^2 = I_E(k) + I_G(k)^2/I_G(k)^2.$$

The map

$$\tilde\varphi \colon R \twoheadrightarrow R_G \hookrightarrow I_G(k)/I_G(k)^2 \xrightarrow{\sim} (G/G^*)\hat\otimes_{\mathbb{F}_p} k \twoheadrightarrow (G/E)\hat\otimes_{\mathbb{F}_p} k$$

is surjective. By Nakayama's lemma it follows that the map

$$\varphi \colon R \hookrightarrow I_G(k) \twoheadrightarrow I_{G/E}(k)$$

is surjective, since $R \xrightarrow{\varphi} I_{G/E}(k) \twoheadrightarrow I_{G/E}(k)/I_{G/E}(k)^2 \cong (G/E)\hat\otimes_{\mathbb{F}_p} k$ is the surjective map $\tilde\varphi$. Thus we have the canonical exact and commutative diagram

$$\begin{array}{ccccccccc}
0 & \longrightarrow & Ind_E^G I_E(k) & \longrightarrow & I_G(k) & \longrightarrow & I_{G/E}(k) & \longrightarrow & 0 \\
& & \uparrow & & \uparrow & & \| & & \\
0 & \longrightarrow & R' & \longrightarrow & R & \xrightarrow{\varphi} & I_{G/E}(k) & \longrightarrow & 0,
\end{array}$$

where $R' = R \cap Ind_E^G I_E(k)$ is the kernel of $\varphi$. Let $\psi \colon I_G(k) \twoheadrightarrow I_{G/E}(k)$ be the map which is equal to $\varphi$ when restricted to $R$ and zero when restricted to $M$. We consider the commutative and exact diagram

$$(*) \qquad \begin{array}{ccccccccc}
& & & & k & = & k & & \\
& & & & \uparrow & & \uparrow & & \\
0 & \longrightarrow & M \oplus R' & \longrightarrow & \Lambda_G(k) & \xrightarrow{\tilde\psi} & X & \longrightarrow & 0 \\
& & \| & & \uparrow & & \uparrow & & \\
0 & \longrightarrow & M \oplus R' & \longrightarrow & I_G(k) & \xrightarrow{\psi} & I_{G/E}(k) & \longrightarrow & 0,
\end{array}$$

where $X$ is defined as the quotient $\Lambda_G(k)/(M \oplus R')$ and $\tilde\psi$ is the projection.

*Claim*: $X \cong \Lambda_{G/E}(k)$.

*Proof*: Let $\{\eta_1, \ldots, \eta_{d-1}\} \subseteq R$ be a pre-image with respect to $\varphi$ of a $k$-basis of $I_{G/E}(k)$. Using (1.8), we get $Res_E I_G(k) = I_E(k) \oplus \bigoplus_{i=1}^{d-1} \Lambda_E(k)\,\eta_i$. Since $\bigoplus_{i=1}^{d-1} \Lambda_E(k)\,\eta_i \subseteq Res_E R$ is a direct summand of $Res_E I_G(k)$, we obtain

$$Res_E R = \tilde R \oplus \bigoplus_{i=1}^{d-1} \Lambda_E(k)\,\eta_i, \quad Res_E R' = \tilde R \oplus \bigoplus_{i=1}^{d-1} I_E(k)\,\eta_i,$$

where $\tilde R = I_E(k) \cap Res_E R$. It follows that

$$Res_E M \oplus \tilde{R} \oplus \bigoplus_{i=1}^{d-1} \Lambda_E(k)\,\eta_i = Res_E M \oplus Res_E R$$
$$= Res_E I_G(k)$$
$$= I_E(k) \oplus \bigoplus_{i=1}^{d-1} \Lambda_E(k)\,\eta_i,$$

and so $Res_E M \oplus \tilde{R} \cong I_E(k)$. Applying the restriction functor to the diagram $(*)$, we obtain the commutative and exact diagram

$$\begin{array}{ccccccccc}
0 & \longrightarrow & Res_E M \oplus \tilde{R} \oplus \bigoplus_{i=1}^{d-1} I_E(k)\,\eta_i & \longrightarrow & \Lambda_E(k) \oplus \bigoplus_{i=1}^{d-1} \Lambda_E(k)\,\eta_i & \longrightarrow & Res_E X & \longrightarrow & 0 \\
& & \| & & \uparrow & & \uparrow & & \\
0 & \longrightarrow & Res_E M \oplus \tilde{R} \oplus \bigoplus_{i=1}^{d-1} I_E(k)\,\eta_i & \longrightarrow & I_E(k) \oplus \bigoplus_{i=1}^{d-1} \Lambda_E(k)\,\eta_i & \longrightarrow & k^{d-1} & \longrightarrow & 0,
\end{array}$$

and, dividing out the term $\bigoplus_{i=1}^{d-1} I_E(k)$, the commutative and exact diagram

$$\begin{array}{ccccccccc}
0 & \longrightarrow & Res_E M \oplus \tilde{R} & \longrightarrow & \Lambda_E(k) \oplus \bigoplus_{i=1}^{d-1} k & \longrightarrow & Res_E X & \longrightarrow & 0 \\
& & \| & & \uparrow & & \uparrow & & \\
0 & \longrightarrow & Res_E M \oplus \tilde{R} & \longrightarrow & I_E(k) \oplus \bigoplus_{i=1}^{d-1} k & \longrightarrow & k^{d-1} & \longrightarrow & 0.
\end{array}$$

Since $Res_E M \oplus \tilde{R} \cong I_E(k)$, we obtain $Res_E X \cong k^d$, and so $Res_E X$ is a trivial $E$-module, i.e. $X$ is a $G/E$-module and we obtain a surjection $\Lambda_{G/E}(k) \twoheadrightarrow X$ which has to be an isomorphism. This proves the claim.

Thus we have two exact sequences of $\Lambda_G(k)$-modules
$$0 \longrightarrow M \oplus R' \longrightarrow \Lambda_G(k) \longrightarrow \Lambda_{G/E}(k) \longrightarrow 0$$
$$0 \longrightarrow Ind_E^G I_E(k) \longrightarrow \Lambda_G(k) \longrightarrow \Lambda_{G/E}(k) \longrightarrow 0.$$
The lemma of Schanuel implies that $M \oplus R' \oplus \Lambda_G(k) \cong Ind_E^G I_E(k) \oplus \Lambda_G(k)$, see [6] (2.24), thus $M | Ind_E^G I_E(k)$ by (1.2). $\square$

**Lemma 1.12** *Let $G$ be a finitely generated pro-$p$ group and $H$ a closed subgroup. Let $M$ be a finitely generated $\Lambda_G$-module such that*
$$M | Ind_H^G I_H.$$
*Assume that $H$ is not finitely generated. Then there exists a proper open subgroup $H_0$ of $H$ such that*
$$M | Ind_{H_0}^G I_{H_0} \quad \text{and} \quad I_{H_0} + I_G^2 = I_H + I_G^2.$$

**Proof:** We may assume that $M$ is contained in $Ind_H^G I_H$. Let $N$ be an open normal subgroup of $G$, $E$ a proper open normal subgroup of $H$ containing $H^*$ and
$$\varphi_N^E \colon M_N \hookrightarrow (Ind_H^G I_H)_N \twoheadrightarrow (Ind_H^G I_{H/E})_N = Ind_{HN/N}^{G/N}(I_{H/E})_{N \cap H}.$$
If $N$ runs through the open normal subgroups of $G$, then the sets
$$\mathcal{B}_N = \{H^* \subseteq E \subsetneq H \mid E \text{ open and normal in } H, \ \varphi_N^E = 0, \ I_E + I_G^2 = I_H + I_G^2\}$$
together with the inclusion maps $\mathcal{B}_{N''} \hookrightarrow \mathcal{B}_{N'}$, $N'' \subseteq N'$, form a projective system. Furthermore, $\mathcal{B}_N = \varprojlim_{H'} \mathcal{B}_{N,H'}$,
$$\mathcal{B}_{N,H'} = \{H' \subseteq E \subsetneq H \mid E \text{ open and normal in } H, \ \varphi_N^E = 0, \ I_E + I_G^2 = I_H + I_G^2\},$$
where $H' \subsetneq H$ is open and normal, $H^* \subseteq H'$, and the transition maps are given by $\mathcal{B}_{N,H''} \to \mathcal{B}_{N,H'}$, $E \mapsto EH'$, if $H'' \subseteq H'$. Since the sets $\mathcal{B}_{N,H'}$ are finite, the set $\mathcal{B}_N$ is compact. We will show that these sets are not empty. Then $\varprojlim_N \mathcal{B}_N$ is not empty, i.e. there exists a proper open normal subgroup $H_0$ of $H$ such that the map
$$\varphi^E \colon M \hookrightarrow Ind_H^G I_H \twoheadrightarrow Ind_H^G I_{H/H_0}$$
is zero. It follows that $M \subseteq Ind_{H_0}^G I_{H_0}$, hence $M | Ind_{H_0}^G I_{H_0}$, and $I_{H_0} + I_G^2 = I_H + I_G^2$.
**Claim:** $\mathcal{B}_N \neq \emptyset$.
Since $M_N$ has finite $\mathbb{F}_p$-dimension, the image of the map
$$M_N \hookrightarrow (Ind_H^G I_H)_N \twoheadrightarrow (Ind_H^G I_{H/H^*})_N = Ind_{HN/N}^{G/N}(I_{H/H^*})_{N \cap H}$$
is contained in $Ind_{HN/N}^{G/N} A$, where $A \subseteq (I_{H/H^*})_{N \cap H}$ is a $H/N \cap H$-module and finite-dimensional as $\mathbb{F}_p$-vector space. Let $\tilde{N} = (N \cap H)H^*$. Consider the commutative and exact diagram

$$\begin{array}{ccccccccc}
0 & \longrightarrow & \tilde{N}/H^* & \longrightarrow & (I_{H/H^*})_{\tilde{N}/H^*} & \longrightarrow & I_{H/\tilde{N}} & \longrightarrow & 0 \\
& & \uparrow & & \uparrow & & \uparrow & & \\
0 & \longrightarrow & A_0 & \longrightarrow & A & \longrightarrow & A/A_0 & \longrightarrow & 0,
\end{array}$$

where $A_0 = A \cap \tilde{N}/H^*$. Let $S \subseteq H/H^*$ be a complement of $H \cap G^*/H^*$; thus $S$ has finite $\mathbb{F}_p$-dimension and is mapped isomorphically onto $HG^*/G^*$ under the canonical map $H/H^* \twoheadrightarrow HG^*/G^*$. Since $\tilde{N}/H^*$ has infinite $\mathbb{F}_p$-dimension and $A_0$ is finite-dimensional, there exists a proper finite-dimensional subspace $0 \neq \tilde{E} \subseteq H/H^*$ such that $\tilde{E} \subseteq \tilde{N}/H^*$ and $\tilde{E} \cap SA_0 = 0$. Let $E$ be a proper open normal subgroup of $H$ containing $H^*$ such that $SA_0 \subseteq E/H^*$ and $\tilde{E} \oplus E/H^* = H/H^*$. It follows that $H/E \cong \tilde{E}$, $\tilde{N}E = H$ and
$$(I_{H/E})_{\tilde{N}} = I_{H/E}/I_{\tilde{N}E/E} I_{H/E} = I_{H/E}/(I_{H/E})^2 = H/E.$$

Thus the map $A \hookrightarrow (I_{H/H^*})_{\tilde N} \twoheadrightarrow (I_{H/E})_{\tilde N}$ is zero, and so $\varphi_N^E = 0$. Furthermore $EG^*/G^* = HG^*/G^*$, i.e. $I_E + I_G^2 = I_H + I_G^2$. This proves the claim and so the lemma. $\square$

**Proposition 1.13** *Let $G$ be a finitely generated pro-p group. Let*
$$I_G = M \oplus R$$
*be a decomposition of the augmentation ideal $I_G$ into left $\Lambda_G$-modules, where $M$ is indecomposable and not isomorphic to $\Lambda_G$. Then there exists a finitely generated closed subgroup $H$ of $G$ such that*
$$M \cong Ind_H^G I_H \quad \text{and} \quad M + I_G^2 = I_H + I_G^2.$$

**Proof:** We may assume that the decomposition of $I_G$ is non-trivial (otherwise take $H = G$). We consider the set
$$\mathcal{M} = \{U \mid U \text{ a closed subgroup of } G \text{ such that } M | Ind_U^G I_U, M + I_G^2 = I_U + I_G^2\},$$
which is partially ordered when ordered by inclusion and non-empty by (1.11). Furthermore, if $\{U_1 \supseteq U_2 \supseteq \ldots\}$ is a chain in $\mathcal{M}$, then $U = \bigcap_i U_i$ is a lower bound in $\mathcal{M}$; indeed, from [10] (3.7) and (4.2) it follows that $M | Ind_U^G Res_U M$. Since $M | I_G$ and so $Res_U M | Res_U I_G$, we obtain, using (1.8),
$$Res_U M \mid (I_U \oplus P_0), \text{ where } P_0 \text{ is } \Lambda_U\text{-free},$$
hence
$$M \mid (Ind_U^G I_U \oplus P), \text{ where } P \text{ is } \Lambda_G\text{-free}.$$
Using (1.2), it follows that $M | Ind_U^G I_U$, since $M \not\cong \Lambda_G$, and so
$$\dim_{\mathbb{F}_p}(M + I_G^2)/I_G^2 \leq \dim_{\mathbb{F}_p}(Ind_U^G I_U + I_G^2)/I_G^2 = \dim_{\mathbb{F}_p}(I_U + I_G^2)/I_G^2.$$
Since $I_U + I_G^2 \subseteq I_{U_i} + I_G^2 = M + I_G^2$, we obtain $M + I_G^2 = I_U + I_G^2$. Now Zorn's lemma implies that $\mathcal{M}$ has a minimal element $H$ and we have $M | Ind_H^G I_H$ and $M + I_G^2 = I_H + I_G^2$.

The group $H$ is finitely generated, since otherwise, by (1.12), there would exist a proper open subgroup $H_0$ of $H$ such that $M | Ind_{H_0}^G I_{H_0}$ and $M + I_G^2 = I_{H_0} + I_G^2$ which contradicts the minimality of $H$.

Furthermore, $I_H(k)$ is indecomposable for every finite extension $k$ of $\mathbb{F}_p$. Indeed, suppose the contrary, i.e. $I_H(k) = A \oplus B$ is a non-trivial decomposition, where $A$ is an indecomposable $\Lambda_H(k)$-module not isomorphic to $\Lambda_H(k)$, and $k \otimes_{\mathbb{F}_p} M | Ind_H^G A$ (use (1.2)). From (1.11) it follows that $A | Ind_E^H I_E(k)$ and $A + I_H(k)^2 = I_E(k) + I_H(k)^2$ where $E$ is a proper open subgroup of $H$. Thus $k \otimes_{\mathbb{F}_p} M | Ind_E^G I_E(k)$. Using (1.9)(ii), we get $M | Ind_E^G I_E$. Furthermore

$$\begin{aligned}
\dim_{\mathbb{F}_p}(I_H + I_G^2)/I_G^2 &= \dim_k(k \otimes_{\mathbb{F}_p} M + I_G(k)^2)/I_G(k)^2 \\
&\leq \dim_k(Ind_H^G A + I_G(k)^2)/I_G(k)^2 \\
&= \dim_k(A + I_G(k)^2)/I_G(k)^2 \\
&= \dim_k(I_E(k) + I_G(k)^2)/I_G(k)^2 \\
&= \dim_{\mathbb{F}_p}(I_E + I_G^2)/I_G^2.
\end{aligned}$$

hence $M + I_G^2 = I_E + I_G^2$. Thus we get a contradiction to the minimality of $H$.

Now Greens indecomposability theorem for pro-$p$ groups says that $Ind_H^G I_H$ is indecomposable, see [10](6.9), thus $M \cong Ind_H^G I_H$. $\square$

**Theorem 1.14** *Let $G$ be a finitely generated pro-$p$ group and let*

$$I_G = M_1 \oplus \cdots \oplus M_s$$

*be a decomposition into indecomposable left $\Lambda_G$-modules $M_i$. Then there exist freely indecomposable closed subgroups $H_i$ of $G$ such that $Ind_{H_i}^G I_{H_i} \cong M_i$, $i = 1, \ldots, s$, and*

$$G = H_1 * \cdots * H_s.$$

*In particular, $G$ is freely indecomposable if and only if $I_G$ is indecomposable.*

**Proof:** Let $r = \#\{i | M_i \cong \Lambda_G\}$, $t = s - r$ and $M = M_1 \oplus \ldots M_t$ such that $M$ has no free summand (after re-numbering). By (1.13) and (1.10), we obtain

$$I_G = M \oplus Ind_F^G I_F, \qquad M \cong Ind_{H_1}^G I_{H_1} \oplus \cdots \oplus Ind_{H_t}^G I_{H_t},$$

where $F$ is a free pro-$p$ subgroup of $G$ of rank $r$ and $H_i$, $i = 1, \ldots, t$, are closed subgroups of $G$ such that $M_i \cong Ind_{H_i}^G I_{H_i}$ and $M_i + I_G^2 = I_{H_i} + I_G^2$. Since

$$\begin{aligned}
I_G/I_G^2 &= (M_1 + I_G^2)/I_G^2 \oplus \cdots \oplus (M_t + I_G^2)/I_G^2 \oplus (I_F + I_G^2)/I_G^2 \\
&= (I_{H_1} + I_G^2)/I_G^2 \oplus \cdots \oplus (I_{H_1} + I_G^2)/I_G^2 \oplus (I_F + I_G^2)/I_G^2,
\end{aligned}$$

and so

$$G/G^* = H_1 G^*/G^* \oplus \cdots \oplus H_t G^*/G^* \oplus F G^*/G^*,$$

the subgroups $H_1, \ldots, H_t, F$ generate $G$. Let $H_i \cong \mathbb{Z}_p$, $i = t+1, \ldots, s$, such that $F = H_{t+1} * \cdots * H_s$. It follows that the canonical map

$$Ind_{H_1}^G I_{H_1} \oplus \cdots \oplus Ind_{H_s}^G I_{H_s} \overset{can}{\twoheadrightarrow} I_G \cong \bigoplus_{i=1}^s Ind_{H_i}^G I_{H_i}$$

induced by the inclusions is surjective. Since an endomorphism of a finitely generated $\Lambda_G$-module is an isomorphism, the map *can* is bijective. Using (1.6), we see that the homomorphism

$$H_1 * \cdots * H_s \xrightarrow{\sim} G$$

induced by the inclusion is bijective.

It remains to show that the groups $H_i$ are freely indecomposable. Suppose one of the $H_i$ decomposes freely, then by (1.6) the augmentation ideal $I_G$ decomposes in more than $s$ non-trivial summands, which is impossible by the Krull-Schmidt-Azumaya theorem. This finishes the proof of the theorem. □

## 2 Free pro-$p$ products

A finitely generated (abstract) group admits a decomposition into a free product of freely indecomposable groups, called its *Grushko decomposition*. The following theorem is a pro-$p$ analog of this result.

**Theorem 2.1**
(i) *Every finitely generated pro-p group $G$ is the free pro-p product of finitely many freely indecomposable closed subgroups $G_i$, $i = 1, ..., s$, i.e.*

$$G = \underset{i=1}{\overset{s}{\text{\Large $*$}}} G_i.$$

(ii) *Let $G = *_{i \in I} G_i$ be the free pro-p product of pro-p groups $G_i$, $i \in I$, and let $H$ be a finitely generated closed subgroup of $G$ which is freely indecomposable. Then $H \cong \mathbb{Z}_p$ or there exist $j \in I$ and $\sigma \in G$ such that $H \subseteq (G_j)^\sigma$.*

(iii) *Let*

$$G = \underset{i=1}{\overset{n}{\text{\Large $*$}}} G_i * F_t = \underset{j=1}{\overset{m}{\text{\Large $*$}}} H_j * F_u$$

*be two decompositions of the finitely generated pro-p group $G$, where the closed subgroups $H_j$ and $G_i$ are freely indecomposable and not isomorphic to $\mathbb{Z}_p$ and $F_t$ and $F_u$ are free pro-p groups of rank $t$ and $u$, respectively. Then $n = m$, $t = u$ and there are elements $\sigma_1, \ldots, \sigma_n$ in $G$ such that (after possibly re-ordering) $G_i = (H_i)^{\sigma_i}$ for $i = 1, \ldots, n$.*

*In particular, the number $s(G) = n + t$ of freely indecomposable factors of $G$ is an invariant of $G$. We call a decomposition as above a* **Grushko decomposition***.*

**Proof:** (i) In contrast to the theory of abstract groups, for pro-$p$ groups it is easy to see that the generator rank $d(H_1 * H_2)$ of a free pro-$p$ product $H_1 * H_2$ is $d(H_1) + d(H_2)$, see [12](3.9.1),(4.1.4), and so assertion (i) is obvious.

(ii) follows from the pro-$p$ analog of Kurosh' subgroup theorem for *finitely generated closed* subgroups of free pro-$p$ products, see [7] Theorem (9.7) or [8] Theorem (4.4) or [11] Proposition (5.2).

Now we obtain (iii) easily: Using (ii), we get

$$G_i \subseteq (H_j)^{\sigma_j} \subseteq (G_k)^{\tau_k},$$

for $\sigma_j, \tau_k \in G$. It follows that $G_i = (G_k)^{\tau_k}$, hence $G_i = (H_j)^{\sigma_j}$, and so $n = m$. Dividing out the normal closure of the subgroup generated by the $H_i$'s, we get an isomorphism $F_t \cong F_u$, hence $t = u$. $\square$

**Lemma 2.2** *Let*
$$1 \longrightarrow H \longrightarrow G \xrightarrow{\psi} G/H \longrightarrow 1$$
*be an exact sequence of pro-$p$ groups, where $G/H \cong \mathbb{F}_p$.*

(i) *Let $G$ be free of rank $r$ with basis $\{\tau_1, \cdots, \tau_r\}$ such that $\psi$ maps $\tau_1$ to a generator of $G/H$ and $\psi(\tau_i) = 1$ for $i \geq 2$. Then $H$ is free of rank $p(r-1)+1$ with basis*
$$\{\tau_1^p, (\tau_i)^{(\tau_1)^k}, k = 0, \ldots, p-1, i = 2, \ldots, r\},$$
*and there is a $\mathbb{Z}_p[G/H]$-isomorphism*
$$H^{ab} \cong \mathbb{Z}_p \oplus \mathbb{Z}_p[G/H]^{r-1}.$$

(ii) *Let*
$$G = \underset{i=1}{\overset{n}{*}} G_i * L,$$
*$n \geq 1$, and each factor $G_i$ of the free pro-$p$ product is mapped surjectively onto $G/H$ and $L$ to $1$. Then there is a decomposition of $H$ as free pro-$p$ product*
$$H = \tilde{H} * F_d, \quad \tilde{H} = \underset{i=1}{\overset{n}{*}} H_i * \underset{j=0}{\overset{p-1}{*}} L^{(\tau_1)^j},$$
*where $H_i = H \cap G_i$, $\tau_i \in G_i \backslash H_i$, $i = 1, \ldots, n$, and $F_d$ is free of rank $(p-1)(n-1)$ with basis $\{(\tau_i \tau_1^{-1})^{(\tau_1)^k}, i \geq 2, k = 1, \ldots, p-1\}$.*

*Furthermore, if $M$ is the normal subgroup of $H$ generated by $\tilde{H}$ and $\bar{F}_d = F_d M/M$, then there is a $\mathbb{Z}_p[G/H]$-isomorphism*
$$(\bar{F}_d)^{ab} \cong (I_{G/H})^{n-1}.$$

**Proof:** (ii) The pro-$p$ analog of Kurosh' subgroup theorem for open subgroups of free pro-$p$ products, see [12] (4.2.1), implies that $H$ has the asserted structure. In order to find a basis of $F_d$, we first consider the case $G_i =< \tau_i > \cong \mathbb{Z}_p$, $i = 1, \ldots, n$, and $L = 1$. Schreier's subgroup theorem of free groups, see [14] (3.6.2)(a) for the profinite case, implies that

$$u_1 = \tau_1^p, \ u_{i,k} = (\tau_i \tau_1^{-1})^{(\tau_1)^k}, \ u_{i,p-1} = (\tau_1)^{p-1} \tau_i = (\tau_i \tau_1^{-1})^{(\tau_1)^{p-1}} \tau_1^{-p},$$

$k = 0, \ldots, p-2$, $i = 2, \ldots, n$, is a basis of $H$. Changing the basis, we get

$$u_1 = \tau_1^p, \ u_{i,k} = (\tau_i \tau_1^{-1})^{(\tau_1)^k}, \ u_i := u_{i,0} \cdots u_{i,p-2} u_{i,p-1} = \tau_i^p,$$

hence $H_i = < \tau_i^p >$, $i \geq 1$, and $\{(\tau_i \tau_1^{-1})^{(\tau_1)^k}, i \geq 2, k = 1, \ldots, p-1\}$ is a basis of $F_d$.

Dividing out the normal subgroup generated by the subgroups $H_i$ (which is also normal in $G$), we are in the case where $G_1 = \cdots = G_n \cong \mathbb{F}_p$ and $L = 1$. In general, we also obtain this case by dividing out the normal subgroup $M$ (which is also normal in $G$), and get so the desired result.

Considering the explicitly given basis of $F_d$, we get the asserted structure of the $\mathbb{Z}_p[G/H]$-module $(\bar{F}_d)^{ab}$, or one can see this as follows. First observe that $\bar{F}_d \cong F_d$ and so $(\bar{F}_d)^{ab}$ is $\mathbb{Z}_p$-free. We have $G/M = \overset{n}{\underset{i=1}{*}} G_i/H_i$ and therefore the exact sequence

$$1 \longrightarrow \bar{F}_d \longrightarrow \overset{n}{\underset{i=1}{*}} G_i/H_i \longrightarrow G/H \longrightarrow 1,$$

inducing the exact sequence

$$0 \longrightarrow (\bar{F}_d^{ab})_{G/H} \longrightarrow \bigoplus_{i=1}^{n} G_i/H_i \longrightarrow G/H \longrightarrow 0.$$

It follows that $(\bar{F}_d^{ab})_{G/H} \cong \mathbb{F}_p^{n-1}$. Since

$$(\bar{F}_d)^{ab} \cong \mathbb{Z}_p[G]^a \oplus (I_G)^b \oplus \mathbb{Z}_p^c,$$

where

$$(p-1)(n-1) = pa + (p-1)b + c,$$

see [6] §30C, we obtain $(\bar{F}_d)^{ab} \cong (I_G)^b$ and $b = n - 1$.

Finally we see that (i) is a special case of (ii): take $n = 1$, $G_1 =< \tau_1 >$ and $L =< \tau_2, \ldots, \tau_r >$. □

The following theorem is in some sense the converse of lemma (2.2)(ii), i.e. the converse of a special case of Kurosh' subgroup theorem; we will use the notion $N_G(L)$ for the normalizer of a closed subgroup $L$ in $G$.

**Theorem 2.3** *Let $G$ be a finitely generated torsion-free pro-$p$ group and let $H$ be an open normal subgroup such that $G/H \cong \mathbb{F}_p$. Then $H$ and $G$ have a Grushko decomposition as a free pro-$p$ product of closed subgroups of the following form:*

$$G = \underset{i=1}{\overset{n}{\ast}} G_i \ast \underset{j=1}{\overset{m}{\ast}} L_j \ast F_\delta,$$

$$H = \underset{i=1}{\overset{n}{\ast}} H_i \ast \underset{j=1}{\overset{m}{\ast}} \underset{\sigma \in G|H}{\ast} (L_j)^\sigma \ast F_d,$$

$n \geq 0$, $m \geq 0$, *where $G_i = N_G(H_i)$ and*

$H_i$ *is freely indecomposable not isomorphic to $\mathbb{Z}_p$ and $G_i/H_i \cong \mathbb{F}_p$,*
$L_j$ *is freely indecomposable and $N_G(L_j) = L_j$,*
$F_d$ *and $F_\delta$ are free of rank $d = (p-1)(n-1)$ and $\delta = 0$, respectively, if $n \geq 1$, and $d = \delta = 1$ otherwise.*

The following corollary would be an immediate consequence of the main theorem of the introduction, but we have to use it in order to prove this theorem.

**Corollary 2.4** *Let $G$ be a finitely generated torsion-free pro-$p$ group and let $H$ be an open subgroup. Then $H$ is freely decomposable if and only if $G$ is freely decomposable.*

**Remarks:** 1. The assumption that $G$ has to be torsion-free is necessary as the following example shows: Let $G = F_2 \times \mathbb{Z}/p\mathbb{Z}$, where $F_2$ is a free pro-$p$ group of rank 2, and $H = F_2$. Then $H$ is freely decomposable, but $G$ not.

2. Also the opposite case is interesting: The group $\mathbb{Z}/2Z \ast \mathbb{Z}/2\mathbb{Z} = \mathbb{Z}_2 \rtimes \mathbb{Z}/2Z$ is decomposable and has an open subgroup which is freely indecomposable. Kurosh' subgroup theorem shows that this is the only pro-$p$ group with this property.

**Proof of (2.4):** The pro-$p$ analog of Kurosh' subgroup theorem for open subgroups of free pro-$p$ products, see [12] (4.2.1), implies the if-part. In order to show the converse we may assume that $H$ is normal in $G$ of index $p$. Now the result follows from (2.3). □

**Proof of (2.3):** Let
$$H = \underset{i=1}{\overset{n}{\ast}} H_i \ast \underset{\lambda=1}{\overset{\mu}{\ast}} K_\lambda \ast F_u,$$

be a decomposition of $H$, where $H_i$ and $K_\lambda$ are freely indecomposable closed subgroups not isomorphic to $\mathbb{Z}_p$, $N_G(H_i) \neq H_i$, $i = 1, \ldots, n$, and $N_G(K_\lambda) = K_\lambda$, $\lambda = 1, \ldots, \mu$, and $F_u$ is free of rank $u$. Since $H$ is finitely generated and $d(H) = \sum_i d(H_i) + \sum_\lambda d(K_\lambda) + u$, the groups $H_i$ and $K_\lambda$ are finitely generated, too. Let $G_i = N_G(H_i)$. As $N_H(H_i) = H_i$, see [14] (9.1.12), we have $H \cap G_i = H_i$, hence a commutative and exact diagram

$$\begin{array}{ccccccccc} 0 & \longrightarrow & H & \longrightarrow & G & \longrightarrow & \mathbb{F}_p & \longrightarrow & 0 \\ & & \uparrow & & \uparrow & & \| & & \\ 0 & \longrightarrow & H_i & \longrightarrow & G_i & \longrightarrow & \mathbb{F}_p & \longrightarrow & 0. \end{array}$$

It follows from Kurosh' subgroup theorem that $G_i$ is freely indecomposable. Furthermore, since $K_\lambda \neq (K_\lambda)^\sigma$ for a representative $\sigma \in G$ of $\bar{\sigma} \in G/H$, $\bar{\sigma} \neq 1$, it follows from (2.1)(ii) that $(K_\lambda)^\sigma$ is a $H$-conjugate of $K_{\lambda'}$ for some $\lambda' \neq \lambda$. Therefore, we can assume that the above decomposition is of the form

$$H = H' * F_u, \text{ where } H' = \mathop{\text{\Large$*$}}_{i=1}^{n} H_i * \mathop{\text{\Large$*$}}_{j=1}^{m'} \mathop{\text{\Large$*$}}_{\sigma \in G|H} (K_j)^\sigma.$$

Observe that the normal closure $M$ of $H'$ in $H$ is also normal in $G$, and so we obtain an exact sequence $1 \to H/M \to G/M \to \mathbb{F}_p \to 1$. Since $H/M \cong F_u$, the group $(H/M)^{ab}$ is $\mathbb{Z}_p$-free and there is a $\mathbb{Z}_p[G/H]$-isomorphism

$$(H/M)^{ab} \cong \mathbb{Z}_p[G/H]^{u_1} \oplus (I_{G/H})^{u_2} \oplus (\mathbb{Z}_p)^{u_3},$$

where $u = pu_1 + (p-1)u_2 + u_3$, see [6] §30C. It follows that $H$ is of the form

$$H = \tilde{H} * F_{(p-1)u_2+u_3},$$

$$\tilde{H} = \mathop{\text{\Large$*$}}_{i=1}^{n} H_i * \mathop{\text{\Large$*$}}_{\sigma \in G|H} (\mathop{\text{\Large$*$}}_{j=1}^{m'} K_j * F_{u_1})^\sigma = \mathop{\text{\Large$*$}}_{i=1}^{n} H_i * \mathop{\text{\Large$*$}}_{\sigma \in G|H} \mathop{\text{\Large$*$}}_{j=1}^{m} (L_j)^\sigma \,;$$

here $L_j$ is a freely indecomposable factor of the form $K_j$ or $\mathbb{Z}_p$ and $m = m' + u_1$. In the following we put $u = (p-1)u_2 + u_3$ and $M$ now denotes the normal closure of $\tilde{H}$ in $H$. If we consider a subgroup $U$ of $G$ modulo $M$, we denote it by $\bar{U}$. We have a $\mathbb{Z}_p[G/H]$-isomorphism

$$(\bar{F}_u)^{ab} \cong (I_{G/H})^{u_2} \oplus (\mathbb{Z}_p)^{u_3}.$$

Let

$$G_0 = \begin{cases} \mathop{\text{\Large$*$}}_{i=1}^{n} G_i * \mathop{\text{\Large$*$}}_{j=1}^{m} L_j, & \text{if } n > 0, \\ \mathop{\text{\Large$*$}}_{j=1}^{m} L_j * \Gamma, & \text{if } n = 0, \end{cases}$$

where $\Gamma = \mathbb{Z}_p$, and let

$$\varphi : G_0 \longrightarrow G$$

be the homomorphism which is the inclusion on each factor $G_i$ or $L_j$ and in the second case a generator of $\Gamma$ is mapped to a pre-image in $G$ of a generator of $G/H$. We obtain a commutative and exact diagram

$$\begin{array}{ccccccccc} 1 & \longrightarrow & H & \longrightarrow & G & \longrightarrow & G/H & \longrightarrow & 1 \\ & & \uparrow & & \uparrow \varphi & & \| & & \\ 1 & \longrightarrow & H_0 & \longrightarrow & G_0 & \longrightarrow & G/H & \longrightarrow & 1, \end{array}$$

where
$$H_0 = \tilde{H} * F_d,$$

with $d = (p-1)(n-1)$, if $n \geq 1$ and $d = 1$ otherwise. In the first case it follows that there is a $\mathbb{Z}_p[G/H]$-isomorphism

$$(\bar{F}_d)^{ab} \cong (I_{G/H})^{n-1},$$

see (2.2)(ii) and observe that $\bar{F}_d = F_d N/N$ where $N$ denotes the normal closure of $\tilde{H}$ in $H_0$.

*Part 1*: Let $\tau \in H$. Then the subgroup $U = <\tau, F_u> \subseteq H$ is free.

*Proof*: Let $U = V * F_t$, $V = *_{k \in I} V_k$, be a free decomposition of $U$, where $F_t$ is free of rank $t$ and the groups $V_k$ are freely indecomposable not isomorphic to $\mathbb{Z}_p$. Since $G$ is torsion-free and so $U$ is, we have $d(V_k) > 1$ if $V_k$ is not trivial. By the Kurosh subgroup theorem for finitely generated closed subgroups it follows that $V \subseteq M$ (recall that $M$ is the normal closure of $\tilde{H}$ in $H$), hence $UM/M = F_t M/M$. Since $UM/M \cong \bar{F}_u$, it follows that $t \geq u$, and since $d(U) \leq u+1$, we obtain that $U = F_t$, where $u \leq t \leq u+1$.

*Part 2*: $\varphi$ is surjective.
*Proof*: We consider the commutative and exact diagram

$$\begin{array}{ccccccc} & & D & \longrightarrow & C & & \\ & & \uparrow & & \uparrow & & \\ (\tilde{H} * F_u)^{ab} & \longrightarrow & G^{ab} & \longrightarrow & G/H & \longrightarrow & 0 \\ \uparrow & & \uparrow \bar{\varphi} & & \| & & \\ (\tilde{H} * F_d)^{ab} & \longrightarrow & G_0^{ab} & \longrightarrow & G/H & \longrightarrow & 0, \end{array}$$

where the map $\bar{\varphi}$ is induced by $\varphi$. The cokernel $D$ is an image of $\bar{F}_u^{ab}$ and $C = D_G$. Hence we have an exact sequence

$$(\bar{F}_d / \bar{F}_d^*)_G \longrightarrow (\bar{F}_u / \bar{F}_u^*)_G \longrightarrow C/C^* \longrightarrow 0.$$

If $\varphi$ is not surjective, then we need additional generators: Let

$$\varphi' \colon G_0 * F_s \twoheadrightarrow G$$

be surjective, where $F_s$ is free of rank $s$, and a basis of $F_s$ is mapped onto generators $x_1, \ldots, x_s$ of $G$ which are pre-images of a basis of $C/C^*$ and contained in $F_u \backslash \varphi(F_d) \subseteq H$. Let $\sigma \in G_0$ be a pre-image of a generator of $G/H$ under the map $G_0 \xrightarrow{\varphi} G \twoheadrightarrow G/H$, which is contained in $G_1 \subseteq G_0$ if $n \geq 1$ and in $\Gamma$ otherwise. We denote $\varphi'(\sigma)$ also by $\sigma$. Consider the following subgroup

$$F_0 = <\sigma> * F_s$$

of $G_0 * F_s$ and the homomorphism $\psi \colon F_0 \hookrightarrow G_0 * F_s \twoheadrightarrow G \twoheadrightarrow G/H$. Using (2.2)(i), we have the exact sequence

$$1 \longrightarrow E_0 \longrightarrow F_0 \xrightarrow{\psi} G/H \longrightarrow 1,$$

where

$$E_0 = <\sigma^p> * E'_0, \qquad E'_0 = \mathop{\huge *}_{i=0}^{p-1} F_s^{\sigma^i},$$

and a $\mathbb{Z}_p[G/H]$-isomorphism

$$(\bar{E}'_0)^{ab} \cong \mathbb{Z}_p[G/H]^s.$$

Furthermore, let $F = \varphi'(F_0)$ and $E = \varphi'(E_0) = F \cap H$. Then

$$E = <\sigma^p, E'>$$

where $E' \subseteq F_u$ is the pre-image of $\varphi'(E'_0)M/M \subseteq \bar{F}_u$. We have the commutative diagrams

$$
\begin{array}{ccccccccc}
1 & \longrightarrow & H & \longrightarrow & G & \longrightarrow & G/H & \longrightarrow & 1 \\
  &                 & \uparrow &           & \uparrow &           & \| & & \\
1 & \longrightarrow & E & \longrightarrow & F & \longrightarrow & G/H & \longrightarrow & 1 \\
  &                 & \uparrow &           & \uparrow\varphi' &     & \| & & \\
1 & \longrightarrow & E_0 & \longrightarrow & F_0 & \longrightarrow & G/H & \longrightarrow & 1 \\
\end{array}
$$

and

$$
\begin{array}{ccccc}
E_0 = <\sigma^p> * E'_0 & \twoheadrightarrow & E = <\sigma^p, E'> & \longrightarrow & EM/M = \bar{E} \\
 & & \cup & & \cup \\
 & & H = \tilde{H} * F_u & \twoheadrightarrow & H/M = \bar{F}_u.
\end{array}
$$

By part 1 the group $E$ is free. Since $F$ as subgroup of $G$ is torsion-free, we obtain from a theorem of Serre [15] that $F$ is free.

The generator rank of $F$ is $d(F) = 1 + s$, as $F/F^* \cong <\bar{\sigma}> \oplus C/C^*$, where $\bar{\sigma} = \sigma F^*$. It follows that $d(E) = ps+1$, hence $E'_0 \xrightarrow{\sim} E'$ and so $(\bar{E}')^{ab} \cong Z_p[G/H]^s$. By definition of $F_s$ and since $C/C^*$ is a direct summand of $\bar{F}_u/\bar{F}_u^*$, it follows that $(\bar{E}')^{ab} \subseteq (\bar{F}_u)^{ab}$ is a direct $\mathbb{Z}_p$-summand. Thus we have an exact sequence of $\mathbb{Z}_p[G/H]$-modules

$$0 \longrightarrow (\bar{E}')^{ab} \longrightarrow (\bar{F}_u)^{ab} \longrightarrow R \longrightarrow 0,$$

where $R = (\bar{F}_u)^{ab}/(\bar{E}')^{ab}$, which splits as a sequence of free $\mathbb{Z}_p$-modules, i.e. $\text{Ext}^1_{\mathbb{Z}_p}(R, (\bar{E}')^{ab}) = 0$. Therefore

$$\text{Ext}^1_{\mathbb{Z}_p[G/H]}(R, (\bar{E}')^{ab}) = H^1(G/H, \text{Hom}_{\mathbb{Z}_p}(R, (\bar{E}')^{ab})) = 0,$$

as $\text{Hom}_{\mathbb{Z}_p}(R, (\bar{E}')^{ab})$ is a cohomological trivial $\mathbb{Z}_p[G/H]$-module. We obtain a $\mathbb{Z}_p[G/H]$-isomorphism

$$R \oplus \mathbb{Z}_p[G/H]^s \cong (\bar{F}_u)^{ab} \cong (I_{G/H})^{u_2} \oplus (\mathbb{Z}_p)^{u_3},$$

hence $s = 0$. This proves part 2.

*Part* 3: Let $n \geq 1$, then $d(G) = d(G_0)$.
*Proof*: We consider the surjection

$$G_0/G_0^* = \bigoplus_{i=1}^n G_i/G_i^* \oplus \bigoplus_{j=1}^m L_j/L_j^* \twoheadrightarrow G/G^*.$$

Suppose that $d(G) < d(G_0)$. Then one summand in $G_0/G_0^*$ can be replaced by a 1-codimensional subspace and the corresponding map remains to be surjective. The only possible replacement is $G_i/G_i^*$ by $H_i G_i^*/G_i^*$, as otherwise the corresponding map $U \to G$, where $U \subseteq G_0$ is a pre-image of the 1-codimensional subspace of $G_0/G_0^*$ obtained by this replacement, would not induce a surjection onto $H$. Thus there would be a surjection (after re-ordering)

$$\varphi' : G_0' = H_1 * \underset{i=2}{\overset{n}{*}} G_i * \underset{j=1}{\overset{m}{*}} L_j \twoheadrightarrow G,$$

and $n$ has to be bigger or equal to 2. Let $\sigma_j \in G_j \backslash H_j$, $j = 1, \ldots, n$. Using the Kurosh subgroup theorem for open subgroups, we obtain the commutative and exact diagram

$$1 \longrightarrow \underset{i=1}{\overset{n}{*}} H_i * \underset{k=0}{\overset{p-1}{*}} (\underset{j=1}{\overset{m}{*}} L_j)^{(\sigma_2)^k} * F_u \longrightarrow G \longrightarrow \mathbb{F}_p \longrightarrow 1$$

$$1 \longrightarrow \underset{k=0}{\overset{p-1}{*}} (H_1)^{(\sigma_2)^k} * \underset{i=2}{\overset{n}{*}} H_i * \underset{k=0}{\overset{p-1}{*}} (\underset{j=1}{\overset{m}{*}} L_j)^{(\sigma_2)^k} * F_{d'} \longrightarrow G_0' \longrightarrow \mathbb{F}_p \longrightarrow 1,$$

where $u = (p-1)u_2 + u_3$, $d' = (p-1)(n-2)$ and $F_{d'}$ is generated by the elements $(\sigma_j \sigma_2^{-1})^{(\sigma_2)^k}$, $j = 3, \ldots, n$, $k = 1, \ldots, p-1$, see (2.2)(ii).

Let $N'$ be the normal closure of $\ast_{i=2}^n G_i \ast \ast_{j=1}^m L_j$ in $G'_0$ and let $M'$ be its image in $G$. From the commutative diagram above it follows that the kernel of $\varphi'$ is contained in the normal closure of the subgroup generated by $(H_1)^{(\sigma_2)^k-1}$, $k = 1, \ldots, p-1$, and $F_{d'}$, in particular in $N'$. Therefore we get a commutative and exact diagram

$$\begin{array}{ccccccccc}
1 & \longrightarrow & M' & \longrightarrow & G & \longrightarrow & H_1 & \longrightarrow & 1 \\
& & \uparrow & & \uparrow{\varphi'} & & \| & & \\
1 & \longrightarrow & N' & \longrightarrow & G'_0 & \rightrightarrows & H_1 & \longrightarrow & 1.
\end{array}$$

It follows that $G$ is the semi-direct product of $H_1$ and $M'$. We can assume that the element $\sigma_1 \in G_1 \setminus H_1$ is chosen such that its pre-image in $G'_0$ lies inside $N'$. Then $(\sigma_1)^p$ is contained in $H_1 \cap M'$, hence $(\sigma_1)^p = 1$. Thus $G$ contains a torsion element, which is a contradiction. This finishes the proof of part 3.

*Part 4*: Let $n \geq 1$, then $\varphi$ is an isomorphism.
*Proof*: We have a commutative and exact diagram

$$\begin{array}{ccccccccc}
1 & \longrightarrow & H & \longrightarrow & G & \longrightarrow & G/H & \longrightarrow & 1 \\
& & \uparrow & & \uparrow{\varphi} & & \| & & \\
1 & \longrightarrow & H_0 & \longrightarrow & G_0 & \longrightarrow & G/H & \longrightarrow & 1, \\
& & \uparrow & & \uparrow & & & & \\
& & K & = & K & & & &
\end{array}$$

where $K = \operatorname{Ker} \varphi$, $H_0 = \tilde{H} \ast F_d$, $d = (p-1)(n-1)$. Recall that $N$ is the normal closure of $\tilde{H}$ in $H_0$ (which is also normal in $G_0$) and $M$ is the normal closure of $\tilde{H}$ in $H$ (which is also normal in $G$). If we consider a subgroup $U$ of $G_0$ modulo $N$, we denote it by $\bar{U}$. We have a $\mathbb{Z}_p[G/H]$-isomorphism

$$(\bar{F}_d)^{ab} \cong (I_{G/H})^{n-1}.$$

Furthermore, since the inflation map $H^2(H, \mathbb{Q}_p/\mathbb{Z}_p) \xrightarrow{\sim} H^2(H_0, \mathbb{Q}_p/\mathbb{Z}_p)$ is an isomorphism, the exact sequence $1 \to K \to H_0 \to H \to 1$ induces the exact sequence

$$0 \to (K^{ab})_H \to (H_0)^{ab} \to H^{ab} \to 0.$$

Since $(N^{ab})_{H_0} \xrightarrow{\bar{\varphi}}_{\sim} (M^{ab})_H$, the commutative and exact diagram

$$\begin{array}{ccc}
& (K^{ab})_H & \\
& \downarrow & \\
0 \longrightarrow (N^{ab})_{H_0} \longrightarrow (H_0)^{ab} \longrightarrow (\bar{F}_d)^{ab} \longrightarrow 0 \\
\downarrow \sim \quad\quad\quad \downarrow \quad\quad\quad \downarrow \\
0 \longrightarrow (M^{ab})_H \longrightarrow H^{ab} \longrightarrow (\bar{F}_u)^{ab} \longrightarrow 0
\end{array}$$

yields the exact sequence
$$0 \to (K^{ab})_H \to (\bar{F}_d)^{ab} \to (\bar{F}_u)^{ab} \to 0,$$

i.e. the exact sequence
$$0 \to (K^{ab})_H \to (I_{G/H})^{n-1} \to (I_{G/H})^{u_2} \oplus \mathbb{Z}_p^{u_3} \to 0.$$

It follows that $u_3$ has to be zero, since $((I_{G/H})^{n-1})_{G/H} = \mathbb{F}_p^{n-1}$ surjects onto $\mathbb{Z}_p^{u_3}$, and $(K^{ab})_H \cong (I_{G/H})^s$, $s = n - 1 - u_2$, since it has no factor isomorphic to $\mathbb{Z}_p$ or $\mathbb{Z}_p[G/H]$ as contained in $(I_{G/H})^{n-1}$. The exact and commutative diagram

$$\begin{array}{ccccccccc}
1 & \longrightarrow & M & \longrightarrow & G & \longrightarrow & G/M & \longrightarrow & 1 \\
& & \uparrow & & \uparrow & & \uparrow & & \\
1 & \longrightarrow & N & \longrightarrow & G_0 & \longrightarrow & G_0/N & \longrightarrow & 1
\end{array}$$

yields the exact and commutative diagram

$$\begin{array}{ccccccccc}
0 & \longrightarrow & MG^*/G^* & \longrightarrow & G/G^* & \longrightarrow & (G/M)/(G/M)^* & \longrightarrow & 0 \\
& & \uparrow & & \uparrow & & \uparrow & & \\
0 & \longrightarrow & NG_0^*/G_0^* & \longrightarrow & G_0/G_0^* & \longrightarrow & (G_0/N)/(G_0/N)^* & \longrightarrow & 0.
\end{array}$$

Using part 3, it follows that $d(G/M) = d(G_0/N) = n$. Now the exact sequence
$$1 \longrightarrow \bar{F}_{(p-1)u_2} \longrightarrow G/M \longrightarrow G/H \longrightarrow 1$$

induces the exact sequence
$$0 \longrightarrow (\bar{F}_{(p-1)u_2})^{ab}_{G/H} \longrightarrow (G/M)^{ab} \longrightarrow G/H \longrightarrow 0,$$

i.e. the exact sequence
$$0 \longrightarrow \mathbb{F}_p^{u_2} \longrightarrow (G/M)^{ab} \longrightarrow \mathbb{F}_p \longrightarrow 0.$$

Since $(G_0/N)^{ab}$ is elementary abelian, the same is true for its homomorphic image $(G/M)^{ab}$. It follows that of $d(G/M) = u_2 + 1$. Hence $u_2 = n - 1$, and so $(K_H^{ab})_{G/H} = 0$. It follows that $K = 1$, i.e. $\varphi$ is bijective. This proves part 4.

*Part 5*: Let $n = 0$, then $\varphi$ is an isomorphism.

*Proof*: Let $\sigma \in G$ be a pre-image of a generator of $G/H$ and $L = \underset{j=1}{\overset{m}{\ast}} L_j$. Consider the commutative diagram

$$\begin{array}{ccccccccc}
1 & \longrightarrow & \underset{i=0}{\overset{p-1}{\ast}} L^{\sigma^i} \ast F_u & \longrightarrow & G & \longrightarrow & G/H & \longrightarrow & 1 \\
& & \uparrow & & \uparrow \varphi & & \| & & \\
1 & \longrightarrow & \underset{i=0}{\overset{p-1}{\ast}} L^{\sigma^i} \ast \Gamma^p & \longrightarrow & L \ast \Gamma & \longrightarrow & G/H & \longrightarrow & 1.
\end{array}$$

Since $\varphi$ is surjective, and so $\Gamma^p \cong \mathbb{Z}_p \twoheadrightarrow (\bar{F}_u)^{ab}$ is a surjection of $\mathbb{Z}_p[G/H]$-modules, we obtain $u_2 = 0$ and $u_3 \leq 1$. Suppose $u_3 = 0$. Then $H = *_{i=0}^{p-1} L^{\sigma^i}$ and the exact sequence

$$0 \longrightarrow H^{ab} \longrightarrow G/[H,H] \longrightarrow G/H \longrightarrow 0$$

splits, because of $H^2(G/H, H^{ab}) = 0$. It follows that $G/G^* \cong G/H \oplus (H/H^*)_{G/H}$, i.e. $d(G) = d(L) + 1$.

Let $G_1$ be a finitely generated torsion-free and freely indecomposable pro-$p$ group not isomorphic to $\mathbb{Z}_p$ having a surjection $G_1 \twoheadrightarrow \mathbb{Z}_p$ (e.g. $G_1 = \mathbb{Z}_p \oplus \mathbb{Z}_p$), and let $G' = G_1 * G$ and $\psi \colon G' \twoheadrightarrow G$ be the homomorphism which is given by the identity on $G$ and the homomorphism $\pi \colon G_1 \twoheadrightarrow \mathbb{Z}_p \to G$ mapping a generator of $\mathbb{Z}_p$ to $\sigma$. Let $H_1 \subseteq G_1$ be the kernel of the surjection $G_1 \xrightarrow{\pi} G \twoheadrightarrow G/H$. Using the Kurosh subgroup theorem, we get an exact and commutative diagram

$$\begin{array}{ccccccccc}
1 & \longrightarrow & H & \longrightarrow & G & \longrightarrow & G/H & \longrightarrow & 1 \\
& & \uparrow & & \psi \uparrow & & \| & & \\
1 & \longrightarrow & H_1 * H * F_{p-1} & \longrightarrow & G_1 * G & \longrightarrow & G/H & \longrightarrow & 1,
\end{array}$$

Considering the group $G'$ instead of $G$, we are in the situation $n \geq 1$ because of $G_1 = N_{G'}(H_1) \neq H_1$, and we can use the result obtained in that case. Hence

$$G' = G_1 * G \cong G_1 * L$$

which is obviously a contradiction as $d(G') = d(G) + d(G_1) = d(L) + 1 + d(G_1)$. Thus $u_3 = 1$ and it follows that the surjection $H_0 \twoheadrightarrow H$ is bijective, i.e. $\varphi$ is an isomorphism. This finishes the proof of the theorem. $\square$

**Proof of Theorem 2:** The equivalence (i)$\Leftrightarrow$(ii) is theorem (1.14).

In order to prove (iii)$\Rightarrow$(iv) let $H$ be an open subgroup of $G$ such that all open subgroups $H' \subseteq H$ are decomposable, and so by (1.6) the augmentation ideals $I_{H'}$ are decomposable. Using (1.4) (i) and (iv), it follows that $h_1(G) = h_1(H') = \infty$.

Assuming (iv), it follows that there is an open subgroup $H$ of $G$ such that $I_H$ has a direct summand isomorphic to $\Lambda_H$, see (1.4)(ii). Suppose that $I_H \cong \Lambda_H$, then by (1.4)(iii) $H \cong \mathbb{Z}_p$ and $h^1(G) = h^1(H) = 1$, a contradiction. It follows that $I_H \cong M \oplus \Lambda_H$, where $M$ is a non-trivial left $\Lambda_H$-module. From theorem (1.14) and (1.10) it follows that $H = H_0 * F_1$ is freely decomposable, $H_0$ a non-trivial closed subgroup of $H$ and $F_1 \cong \mathbb{Z}_p$. Using Kurosh' subgroup theorem we see that all open subgroups of $H$ are freely decomposable.

If $G$ is torsion free, then by (2.4) we have (i)$\Leftrightarrow$(iii). $\square$

**Proof of Corollary 1:** From theorem 2 and (1.3)(ii) it follows that $G$ is freely indecomposable if and only if $H$ is. Let

$$G = \mathop{\text{\Large *}}_{i=1}^{s(G)} G_i$$

be a decomposition of $G$ into freely indecomposable closed subgroups $G_i$. Recall the pro-$p$ analog of Kurosh' subgroup theorem for open subgroups of free pro-$p$ products, see [12] (4.2.1): There exist systems $S_i$ of representatives $s_i$ of the double coset decomposition $G = \bigcup_{s_i \in S_i} H s_i G_i$ for all $i$ and a free pro-$p$ group $F_r$ of the finite rank

$$r = \sum_{i=1}^{s(G)} [(G:H) - \#S_i] - (G:H) + 1,$$

such that the natural inclusions induce a free product decomposition

$$H = \mathop{\text{\Large *}}_{i=1}^{s(G)} \mathop{\text{\Large *}}_{s_i \in S_i} (G_i^{s_i} \cap H) * F_r.$$

Since $G_i^{s_i} \cap H$ is an open subgroup of $G_i^{s_i}$, it is freely indecomposable by (2.4), hence

$$s(H) = \sum_{i=1}^{s(G)} \#S_i + \sum_{i=1}^{s(G)} [(G:H) - \#S_i] - (G:H) + 1 = (G:H)(s(G) - 1) + 1.$$

This proves part (ii) of the corollary. The proof of part (i) follows from proposition (2.5) below. □

**Proof of Corollary 2:** Recall that $G$ is a duality group of dimension 2 if and only if $D_1(G, \mathbb{F}_p) = 0$, i.e. $h^1(G) = 0$, see [12] (3.4.6). According to theorem 2 this is equivalent to the indecomposability of $G$ resp. $I_G$.

Now let $G$ be a torsion-free 2-generator group. If it would be freely decomposable, then it would be free. Thus it is freely indecomposable if $\text{cd}_p G = 2$, and so it is a duality group.

**Proposition 2.5** *Let $G$ be a finitely generated torsion-free pro-p group. Let*

$$G = \mathop{\text{\Large *}}_{i=1}^{t} H_i * F_r.$$

*be a decomposition of $G$ as free pro-p product, where $H_i$, $i = 1, \ldots, t$, are freely indecomposable closed subgroups, which are not free, and $F_r$ is a free group of rank $r$, i.e. $s(G) = t + r$. Let*

$$I_G = \bigoplus_{i=1}^{t'} M_i \oplus P_{r'}$$

be a decomposition of $I_G$ into indecomposable left $\Lambda_G$-modules $M_i$, $i = 1, \ldots, t'$, which are not isomorphic to $\Lambda_G$, and a free $\Lambda_G$-module $P_{r'}$ of rank $r'$. Then we have the following assertions.

(i) Let $G$ be free. Then $t = t' = 0$, $r = r' = d(G)$,
$$s(G) = f(G),$$
the $\Lambda_G$-module $\mathrm{Hom}_{\Lambda_G}(I_G, \Lambda_G)$ is free of rank $d(G) = f(G)$ and there is an exact sequence $0 \longrightarrow \Lambda_G \longrightarrow \Lambda_G^{f(G)} \longrightarrow H^1(G, \Lambda_G) \longrightarrow 0$.

(ii) If $G$ is not free, then $t = t'$, $r = r'$,
$$s(G) = f(G) + 1$$
and the $\Lambda_G$-modules $\mathrm{Hom}_{\Lambda_G}(I_G, \Lambda_G)$ and $H^1(G, \Lambda_G)$ are free of rank $f(G)+1$ and $f(G)$, respectively. Furthermore, up to a re-ordering there are isomorphisms $M_i \cong J_{H_i}$, $i = 1, \ldots, t$, of left $\Lambda_G$-modules.

**Proof:** Since $I_{F_r} \cong (\Lambda_{F_r})^r$, [12] (5.6.3), (5.6.4), we have $J_{F_r} \cong (\Lambda_G)^r$. Using (1.6), we obtain
$$I_G = \bigoplus_{i=1}^{t} J_{H_i} \oplus (\Lambda_G)^r.$$

Since $H_i$ is not free, we have
$$\mathrm{Ext}^1_G(J_{H_i}, \mathbb{F}_p) \cong \mathrm{Ext}^1_{H_i}(I_{H_i}, \mathbb{F}_p) \cong H^2(H_i, \mathbb{F}_p) \neq 0,$$
hence $J_{H_i}$ is not a free $\Lambda_G$-module.

If $G$ is free of rank $d(G)$, then we have an isomorphism $I_G \cong \Lambda_G^{d(G)}$, [12] (5.6.3), (5.6.4). Furthermore, since $G$ is a duality group with dualizing module $\varinjlim_n D_1(G, \mathbb{Z}/p^n\mathbb{Z}))$ and since $H^1(G, \Lambda_G) = D_1(G, \mathbb{Z}/p\mathbb{Z})^\vee$, we obtain
$$H^1(G, \Lambda_G)_G = H^0(G, D_1(G, \mathbb{Z}/p\mathbb{Z}))^\vee \cong H^1(G, \mathbb{F}_p).$$

If follows that $d(G) = \dim_{\mathbb{F}_p} H^1(G, \Lambda_G)_G = f(G)$. This proves (i).

If $G$ is not free, then $t \geq 1$. From theorem 2 it follows that $h^1(H_i) = 0$, $i = 1, \ldots, t$. Since $\Lambda_G$ is $\Lambda_{H_i}$-projective, the functor $\Lambda_G \hat{\otimes}_{\Lambda_{H_i}} -$ is exact, and we obtain the exact sequence

(†) $$0 \longrightarrow J_{H_i} \longrightarrow \Lambda_G \longrightarrow \Lambda_G \hat{\otimes}_{\Lambda_{H_i}} \mathbb{F}_p \longrightarrow 0.$$

Using (1.5)(ii) and the assumption that $G$ is torsion-free and so $H_i$ is not finite, we have
$$\mathrm{Hom}_{\Lambda_G}(\Lambda_G \hat{\otimes}_{\Lambda_{H_i}} \mathbb{F}_p, \Lambda_G) \cong \mathrm{Hom}_{\Lambda_{H_i}}(\mathbb{F}_p, \mathrm{Res}_{H_i} \Lambda_G) = 0$$
and
$$\mathrm{Ext}^1_{\Lambda_G}(\Lambda_G \hat{\otimes}_{\Lambda_{H_i}} \mathbb{F}_p, \Lambda_G) \cong \mathrm{Ext}^1_{\Lambda_{H_i}}(\mathbb{F}_p, \mathrm{Res}_{H_i} \Lambda_G) = H^1(H_i, \mathrm{Res}_{H_i} \Lambda_G) = 0.$$

Therefore the exact sequence (†) yields the isomorphism

$$\Lambda_G \xrightarrow[\sim]{\tilde{\phi}_{H_i}} \mathrm{Hom}_{\Lambda_G}(J_{H_i}, \Lambda_G),$$

where $\tilde{\phi}_{H_i} = \Lambda_G \hat{\otimes}_{\Lambda_{H_i}} \phi_{H_i}$. Thus $J_{H_i}$ is indecomposable and $\mathrm{Hom}_{\Lambda_G}(I_G, \Lambda_G) \cong (\Lambda_G)^{t+r}$. Furthermore, it follows that the composite map

$$\Lambda_G \xrightarrow{\phi_G} \mathrm{Hom}_{\Lambda_G}(I_G, \Lambda_G) \xrightarrow{pr_1} \mathrm{Hom}_{\Lambda_G}(J_{H_1}, \Lambda_G)$$

is an isomorphism. Therefore, we get an isomorphism

$$(\Lambda_G)^{t-1+r} \cong \bigoplus_{i=2}^{t} \mathrm{Hom}_{\Lambda_G}(J_{H_i}, \Lambda_G) \oplus \mathrm{Hom}_{\Lambda_G}((\Lambda_G)^r, \Lambda_G) \xrightarrow{\sim} H^1(G, \Lambda_G).$$

The last assertion follows from the Krull-Schmidt-Azumaya theorem (1.1)(ii).
□

## 3 Appendix: Abstract groups

A finitely generated (abstract) group admits a decomposition into a free product of freely indecomposable groups, called its *Grushko decomposition*:

*Theorem (Grushko, Kurosh): Every finitely generated group $G$ is the free product of finitely many freely indecomposable subgroups $G_i$, $i = 1, ..., s$, i.e.*

$$G = \overset{s}{\underset{i=1}{*}} G_i.$$

*This decomposition is unique in the following sense: if*

$$G = \overset{n}{\underset{i=1}{*}} G_i = \overset{m}{\underset{j=1}{*}} H_j$$

*are two decompositions of $G$ into freely indecomposable subgroups, then $n = m$ and if $G_i$ is not isomorphic to $\mathbb{Z}$, then there exist an element $\sigma_i$ in $G$ such that (after possibly re-ordering)*

$$G_i = (H_i)^{\sigma_i}.$$

*In particular, the number $s(G) = n$ of freely indecomposable factors of $G$ is an invariant of $G$.*

The uniqueness statement above follows from Kurosh' subgroup theorem.

*Theorem (Kurosh): Let*

$$G = \overset{s}{\underset{i=1}{*}} G_i.$$

*be a decomposition of the finitely generated group $G$ and let $H$ be a subgroup of $G$ of finite index. Then $H$ admits a free product decomposition*

$$H = \overset{s}{\underset{i=1}{*}} \underset{s_i \in S_i}{*} (G_i^{s_i} \cap H) * F_r,$$

*where $S_i$ are systems of representatives $s_i$ of the double coset decomposition $G = \bigcup_{s_i \in S_i} H s_i G_i$ and $F_r$ is a free group of the finite rank*

$$r = \sum_{i=1}^{s} [(G:H) - \#S_i] - (G:H) + 1.$$

Let $R$ be an arbitrary principal ideal domain.

*Theorem: Let $G$ be a finitely generated torsion-free group.*

*(i) (Hopf): The $R$-module $H^1(G, R[G])$ is a free of rank $0, 1$ or $\infty$ and $\mathrm{rank}_R H^1(G, R[G]) = 1$ if and only if $G = \mathbb{Z}$.*

*(ii) (Stallings): $\mathrm{rank}_R H^1(G, R[G]) = \infty$ if and only if $G$ is freely decomposable.*

Since $H^1(G, R[G]) \cong H^1(H, R[H])$ if $H$ is a subgroup of $G$ of finite index, we get the following corollary.

*Corollary: Let $G$ be a finitely generated torsion-free group and let $H$ be a subgroup of finite index. Then $H$ is freely decomposable if and only if $G$ is freely decomposable.*

The following theorem is an immediate consequence of the theorems above and might be well-known but we can not find it in the literature.

**Theorem 3.1** *Let $G$ be a finitely generated torsion-free group and let $H$ be a subgroup of finite index. Then*

$$s(H) = (G:H)(s(G) - 1) + 1,$$

*where $s(G)$ and $s(H)$ are the number of freely indecomposable factors of $G$ and $H$, respectively.*

**Proof:** Let $G = *_{i=1}^{s(G)} G_i$ be a decomposition of $G$ into freely indecomposable subgroups $G_i$. By Kurosh' subgroup theorem we have the decomposition

$$H = \overset{s(G)}{\underset{i=1}{*}} \underset{s_i \in S_i}{*} (G_i^{s_i} \cap H) * F_r.$$

Since $G_i^{s_i} \cap H$ is a subgroup of $G_i^{s_i}$ of finite index, it is freely indecomposable, hence

$$s(H) = \sum_{i=1}^{s(G)} \#S_i + \sum_{i=1}^{s(G)}[(G:H) - \#S_i] - (G:H) + 1 = (G:H)(s(G) - 1) + 1.$$

$\square$

Concerning the structure of $H^1(G, R[G])$ one knows that this $R$-module is isomorphic to $0$ or $R$ or $\bigoplus_1^\infty R$, if $G$ is an infinite finitely generated group. Now we consider this cohomology group with its right $R[G]$-module structure given by right multiplication on $R[G]$. Again we define

$$f(G) = \mathrm{rank}_R \; H^1(G, R[G])_G$$

(we will see that $f(G)$ does not depend on $R$).

**Theorem 3.2** *Let $G$ be a finitely generated torsion-free group. Let*

$$G = \underset{i=1}{\overset{s(G)}{*}} G_i.$$

*be a decomposition of $G$ as free product, where $G_i$, $i = 1, \ldots, s(G)$, are freely indecomposable subgroups. Then we have the following assertions.*

(i) *If $G$ is free, then $s(G) = f(G)$, and there is an exact sequence*

$$0 \longrightarrow R[G] \longrightarrow R[G]^{f(G)} \longrightarrow H^1(G, R[G]) \longrightarrow 0.$$

(ii) *If $G$ is not free, then $s(G) = f(G) + 1$ and the right $R[G]$-module $H^1(G, R[G])$ is free of rank $f(G)$.*

**Proof:** Since $G$ is infinite, the exact sequence

$$0 \longrightarrow I_G \longrightarrow R[G] \longrightarrow R \longrightarrow 0$$

yields the exact sequence

$$0 \longrightarrow R[G] \overset{\phi}{\longrightarrow} \mathrm{Hom}_{R[G]}(I_G, R[G]) \longrightarrow H^1(G, R[G]) \longrightarrow 0,$$

where $\phi(\xi) : x \mapsto \xi x$. Using [5] (4.7), we have

$$I_G = \bigoplus_{i=1}^{s(G)} J_{G_i},$$

where $J_{G_i} = R[G] \otimes_{R[G_i]} I_{G_i}$.

If $G$ is free of rank $s(G)$, then $I_G \cong R[G]^{s(G)}$, as $I_{\mathbb{Z}} \cong R[\mathbb{Z}]$, and we obtain the exact sequence

$$0 \longrightarrow R[G] \longrightarrow R[G]^{s(G)} \longrightarrow H^1(G, R[G]) \longrightarrow 0.$$

Since $G$ is a duality group with dualizing module $H^1(G, \mathbb{Z}[G])$, see [2] (5.1), we obtain

$$H_0(G, H^1(G, R[G])) \cong H^1(G, R) \cong R^{s(G)}.$$

If follows that $s(G) = f(G)$.

If $G$ is not free, then at least one factor $G_i$ is not isomorphic to $\mathbb{Z}$, say $G_1$. From Stallings' theorem it follows that $H^1(G_1, R[G_1])) = 0$. Since $R[G]$ is $R[G_1]$-projective, the functor $R[G] \otimes_{R[G_1]} -$ is exact, and we obtain the exact sequence

(†) $$0 \longrightarrow J_{G_1} \longrightarrow R[G] \longrightarrow R[G] \otimes_{R[G_1]} R \longrightarrow 0.$$

Using Frobenius reciprocity, we have

$$\mathrm{Hom}_{R[G]}(R[G] \otimes_{R[G_1]} R, R[G]) \cong \mathrm{Hom}_{R[G_1]}(R, \mathrm{Res}_{G_1} R[G]) = 0$$

and

$$\mathrm{Ext}^1_{R[G]}(R[G] \otimes_{R[G_1]} R, R[G]) \cong \mathrm{Ext}^1_{R[G_1]}(R, \mathrm{Res}_{G_1} R[G]) = H^1(G_1, \mathrm{Res}_{G_1} R[G]) = 0.$$

From the exact sequence (†) we obtain the isomorphism

$$R[G] \xrightarrow[\sim]{\phi} \mathrm{Hom}_{R[G]}(J_{G_1}, R[G]),$$

and so a commutative and exact diagram

$$\begin{array}{ccc}
\mathrm{Hom}_{R[G]}(\bigoplus_{i=2}^{s(G)} J_{G_i}, R[G]) & \longrightarrow & H^1(G, R[G]) \\
\downarrow & & \| \\
0 \longrightarrow R[G] \longrightarrow \mathrm{Hom}_{R[G]}(I_G, R[G]) & \longrightarrow & H^1(G, R[G]) \longrightarrow 0 \\
\| & \downarrow & \\
R[G] \xrightarrow{\sim} \mathrm{Hom}_{R[G]}(J_{G_1}, R[G]). & &
\end{array}$$

Since $\mathrm{Hom}_{R[G]}(J_{G_i}, R[G]) \cong R[G]$ for all factors $G_i$ (obviously for free factors and for non-free be the consideration above), we get the isomorphism

$$H^1(G, R[G]) \cong \bigoplus_{i=2}^{s(G)} \mathrm{Hom}_{R[G]}(J_{G_i}, R[G]) \cong R[G]^{s(G)-1}.$$

This finishes the proof of the theorem. $\square$

Since in the situation of abstract infinite groups the Krull-Schmidt-Azumaya theorem does not hold in general, and furthermore the $R[G]$-module $R[G]$ is not necessarily indecomposable, we have not the full analog to the assertion (2.5).

# References

[1] Anderson, F.W., Fuller, K.R. *Rings and Categories of Modules.* 2nd edition, Springer 1992

[2] Bieri, R., Eckmann, B. *Groups with Homological Duality Gerneralizing Poincaré Duality.* Invent. math. **20** (1973) 103-124

[3] Brumer, A. *Pseudocompact algebras, profinite groups and class formations.* J. of Algebra **4** (1966) 442–470

[4] Cartan, E., Eilenberg, S. *Homological Algebra.* Princeton Math. Ser. **19**, Princeton 1956

[5] Cohen, D. E. *Groups of Comological Dimension One.* Lect. Notes in Math. **245**, Springer Berlin-Heidelberg-New York 1972

[6] Curtis, C. W., Reiner, I. *Methods of Representation Theory Vol I.* Wiley-Interscience Publication New York, Chichester, Brisbane, Toronto 1981

[7] Haran, D. *On closed subgroups of free products of profinite groups.* Proc. Lond. Math. Soc. **55** (1987), 266-298

[8] Herfort, W.N., Ribes, L. *Subgroups of free pro-p products.* Math. Proc. Camb. Phil. Soc. **101** (1987), 197-206

[9] Korenev, A. A. *Pro-p Groups with Finite Number of Ends.* Math. Notes **76** (2004) 490–496

[10] MacQuarrie, J. *Modular representations of profinite groups* J. pure and applied algeba **215** (2011) 753–763

[11] Melnikov, O.V. *Subgroups and homologies of free products of profinite groups.* Math. USSR Izv. **34** (1990) 97-119

[12] Neukirch, J., Schmidt, A., Wingberg, K. *Cohomology of Number Fields.* 2nd edition, Springer 2008

[13] Pletch, A. *Profinite duality groups II.* J. Pure Applied Algebra **16** 1980, 285–297

[14] Ribes, L., Zalesskiĭ, P. A. *Profinite Groups.* 2nd edition, Springer 2010

[15] Serre, J.-P. *Sur la dimension cohomologique des groupes profinis.* Topology **3** (1965) 413–420

[16] Stallings, J. *On torsion-free groups with infinitely many ends.* Ann. of Math. **88** (1968) 312–334

[17] Swan, R., G. *Groups of Cohomological Dimension One.* J. of Algebra **12** (1969) 585–601

[18] Symonds, P. *On the construction of permutation complexes for profinite groups.* Geometry & Topology Monographs **11** (2007) 369–378

[19] Weigel, Th., Zalesskiĭ, P. A. *Stallings' decomposition theorem for finitely generated pro-p groups.* Preprint 2013

Mathematisches Institut
der Universität Heidelberg
Im Neuenheimer Feld 288
69120 Heidelberg
Germany

e-mail: wingberg at mathi dot uni-heidelberg dot de